\documentclass[11pt,reqno]{amsart} 
\usepackage{amsmath, amssymb, amsthm, mathrsfs, hyperref, a4wide, xcolor}
\usepackage[normalem]{ulem}
\usepackage{url}
\usepackage{comment}

\usepackage{algorithm}
\usepackage{algorithmic}
\usepackage{float}
\usepackage{tikz}
\usetikzlibrary{positioning, arrows.meta, shapes.geometric, fit, backgrounds, calc}
\usepackage{cleveref}
\usepackage[backend=biber,style=numeric]{biblatex}
\DeclareFieldFormat[article,inproceedings,incollection]{title}{#1}
\DeclareFieldFormat[article]{pages}{#1}
\renewbibmacro{in:}{%
  \ifentrytype{article}
    {}
    {\printtext{\bibstring{in}\addspace}}%
}

\newcommand{\EE}{\mathbb{E}}

\newcommand{\dlim}{\xrightarrow[n\to\infty]{d}}

\newcommand{\lp}{\left(}
\newcommand{\Var}{\mathbb{V}}
\newcommand{\rp}{\right)}

\newcommand{\Ord}[1]{{\mathcal O}\left(#1\right)}
\newcommand{\rcp}[1]{\frac{1}{#1}}

\newcommand{\cE}{{E}}
\newcommand{\bt}{\begin{theorem}}
\newcommand{\et}{\end{theorem}}
\newcommand{\bl}{\begin{lemma}}
\newcommand{\el}{\end{lemma}}
\newcommand{\bp}{\begin{prop}}
\newcommand{\ep}{\end{prop}}
\newcommand{\bdf}{\begin{df}}
\newcommand{\edf}{\end{df}}
\newcommand{\bpf}{\begin{proof}}
\newcommand{\epf}{\end{proof}}
\newcommand{\eps}{\varepsilon}
\newcommand{\N}{\mathbb N}
\newcommand{\nti}{n\to\infty}
  
\newcommand{\LambertW}{W} 
\newcommand{\sech}{\operatorname{sech}}
\newcommand{\csch}{\operatorname{csch}}

\newcommand{\ie}{\emph{i.e.}, }

\newcommand{\fra}[1]{\textcolor{blue}{#1}}

\newtheorem{theorem}{Theorem}
\newtheorem{lemma}{Lemma}

\newtheorem{prop}{Proposition}

\theoremstyle{remark}
\newtheorem{remark}{Remark}
\newtheorem{df}{Definition}

\title[Number of histories in two random graph evolution models]
{Asymptotic enumeration of histories and limit laws in two random graph
evolution models}
\thanks{This work has been supported by the Austrian Science Foundation FWF, grant I-6744-N and the Agence Nationale de Recherche, grant ANR-23-CE48-0014. \\
A preliminary version of this work was published in~\cite{AzzouzBodiniDurandGittenberger2026}.}

\author{Nadja Azzouz}
\address{Institute of Discrete Mathematics and Geometry, TU Wien, Austria}
\email{nadja.azzouz@tuwien.ac.at}
\thanks{ORCID Azzouz: \href{https://orcid.org/0009-0004-6060-3476}{0009-0004-6060-3476}}

\author{Olivier Bodini}
\address{Institut Galilée, Univ. Paris Nord}
\email{olivier.bodini@univ-paris13.fr}
\thanks{ORCID Bodini: \href{https://orcid.org/0000-0002-1867-667X}{0000-0002-1867-667X}}

\author{Francis Durand}
\address{Institut Galilée, Univ. Paris Nord}
\email{francis.durand@ens-paris-saclay.fr}
\thanks{ORCID Durand: \href{https://orcid.org/0009-0004-4146-0289}{0009-0004-4146-0289}}

\author[Bernhard Gittenberger]{Bernhard Gittenberger$^*$}
\address{Institute of Discrete Mathematics and Geometry, TU Wien, Austria} 
\email{gittenberger@dmg.tuwien.ac.at}
\thanks{ORCID Gittenberger: \href{https://orcid.org/0000-0002-2639-8227}{0000-0002-2639-8227}}
\thanks{$^*$ Corresponding author}

\date{\today}

\begin{document}

\begin{abstract}
We study growth models for random graphs where the growth is established by vertex or edge creations and where there is no dependence between vertices and edges. Creations are logged by time stamps such that the history of the resulting graph can be reconstructed.  We investigate two prototypical models, an undirected and a directed one, called undirected independent growth model (UIG-model) and directed independent growth model (DIG-model). The generating functions for the number of histories have radius of convergence equal to zero. In the UIG-model, we express the asymptotic number of histories as a sum that is amenable to the discrete Laplace method. We derive the asymptotic number of histories, study the number of vertices of a random history, and prove a central limit theorem. 

In the DIG-model we use a rescaled Borel transform to get an entire function. Then we show that this function is Hayman admissible, which enables us to derive the asymptotic number of histories in this model. The counting sequence is closely related to the framework of generalized Legendre-Stirling numbers. To determine the number of vertices in a random history, we prove that the bivariate generating function is e-admissible and infer that again a normal limit law holds. Finally, we discuss a correspondence to set partitions with even block sizes and present an algorithm based on Boltzmann samplers that generates such graphs uniformly at random. 
\end{abstract}

\maketitle

\section{Introduction}
\label{intro}
Dynamic graph models have emerged as useful tools for modeling systems whose interaction structure changes over time. Examples include social networks, where the times at which interactions occur form an essential part of the data~\cite{berger2006framework}, and stochastic graph-rewriting models for biochemical systems~\cite{behr2020rewriting}. More broadly, random graph theory contains a wide variety of static and evolving models. It started with the classical Erdős--Rényi random graph process~\cite{erdHos1960evolution}, which has been extensively studied and is nowadays well-covered by~\cite{Bollobas01, JLR00} and ranges to more recent models for social networks like the preferential attachment model surveyed by Albert and Barabási~\cite{albert2002statistical}.

In the classical dynamic models there is typically no or no strict separation between the growth in the number of vertices and the growth in the number of edges, see~\cite{albert2002statistical, CooperFrieze03}. For instance, in the Erdős--Rényi models the number of vertices, say $n$, is fixed and the number of edges grows. Then one lets $n$ tend to infinity and studies asymptotic properties, where a dependence on $n$ for the number of edges is usually imposed on the model. In~\cite{farber2024random} a new random graph growth model is suggested where the size as well as the order of a graph may grow according to the inherent structure of the graph process itself. Classical accounts on dynamic random graphs focus on the study of several popular models such as the copying model, the small-world graph, and the preferential attachment model as well as variations in~\cite{Bonato08}, power law graphs and duplication models (static and dynamic) in~\cite{ChungLu06} and the treatment largely based branching processes in~\cite{Durrett07}. 

In this paper, we study a process on multi-graphs which may grow \emph{independently} with respect to order as well as size and where the vertices and the edges receive a time stamp. The evolution is accomplished through one of two operations at each discrete time step: the creation of a new vertex or the addition of an edge between two existing vertices. There are several ways to do that and in this work we consider two particular models: the \emph{undirected independent growth model} (UIG-model) with undirected edges and where no loops are allowed, and the \emph{directed independent growth model} (DIG-model) with directed edges and where loops are allowed. 

For both models, the counting sequence is related to the Legendre–Stirling numbers introduced in~\cite{everitt2002legendre}. The counting sequence for the UIG-model is, up to a one-step shift, OEIS A124373, while the DIG-model gives the sequence of generalized Bell numbers, OEIS A135920. These sequences belong to a broader algebraic framework of polynomial base changes that traces back to Carlitz~\cite{carlitz1963divided} and arises naturally in spectral analysis~\cite{everitt2002legendre}. It was later systematically extended to the Legendre–Stirling framework by Andrews and Littlejohn~\cite{andrews2009combinatorial,andrews2011legendre}, who also gave a combinatorial interpretation in terms of certain set partitions over sets of integers present in two colors. Across these families, the counting sequences satisfy analogous recurrence relations, and the functional equations for their generating functions differ only by a single parameter. Furthermore, they all admit natural graph-theoretic interpretations as special instances of the diagram families studied in~\cite{blasiak2010combinatorial}.

Our contributions are twofold. First, we derive explicit formulas and asymptotic expressions for the number of random graphs obtainable in $n$ steps. 
We regard these graphs as combinatorial structures in the sense of analytic combinatorics and use the machinery developed in \cite{flajolet2009analytic} to derive differential equations for the generating functions associated with these graphs, which are then analyzed asymptotically. The generating functions we encounter have a radius of convergence equal to zero. In the UIG-model, we express the coefficients as sums involving Bessel functions. The terms of the sum are Gaussian and therefore the asymptotic analysis can be based on the discrete Laplace method (see~\cite[p.~761--762]{flajolet2009analytic}). For the DIG-model, we transform the generating function with a rescaled Borel transform~\cite{Balser1994} (see also~\cite{Ovidiu2008} for a general introduction to Borel transforms), making it entire after all. For the transformed function we obtain an explicit closed form. We show that it is Hayman admissible~\cite{Hayman1956} in one variable.

Second, we analyze the distribution of the number of vertices that are created after $n$ steps in both models. In the UIG-model, a perturbed discrete Laplace method applies and so we show a central limit theorem. The Borel transform of the bivariate function of the DIG-model is shown to be e-admissible~\cite{DrmotaGittenbergerKlausner2005}, which yields a normal limit law. We note that in the extended abstract of this work~\cite{AzzouzBodiniDurandGittenberger2026}, a custom modification of bivariate Hayman admissibility~\cite{GittenbergerMandlburger2006} was developed to prove the limit theorem. Bivariate Hayman admissibility is an inherently stronger concept, as it establishes a full bivariate normal distribution. However, since our generating function narrowly misses its conditions, we opt for e-admissibility in this work, which presents a cleaner and more direct approach to derive the result without the need for a conceptual extension.  

\medskip
The paper is organized as follows. Section~\ref{modelA} presents the UIG-model and is devoted to the asymptotic number of graphs that can be generated in $n$ growth steps. Section~\ref{sec:params-undirected} is devoted to the central limit theorem for the distribution of the number of vertices in that model. Out of a desire to tune the growth rate of the number of vertices relative to that for the number of edges, we briefly discuss a weighted version of this the UIG-model in Section~\ref{sec:weighted-undirected}. The DIG-model is introduced in Section~\ref{sec:directed}, where we derive the asymptotic number of graphs that can be generated in $n$ steps. The central limit theorem for the number of vertices is shown in Section~\ref{param_B}. An interesting correspondence to set partitions into even-sized subsets is presented in Section~\ref{subsec:directed-even-partitions}. 
This correspondence can be used to design a uniform sampling algorithm based on Boltzmann samplers~\cite{Duchon2004}, which we present in Section~\ref{algorithm}.

\subsubsection*{Notations and facts}
Throughout the paper we write $\partial_x$ for denoting the partial derivative with respect to $x$. For sequences $(a_n)_{n\in\N}$ and $(b_n)_{n\in\N}$ we write $a_n\asymp b_n$ synonymous to $a_n=\Theta(b_n)$. 

We denote by $j_{\nu,1}$ the smallest positive zero of the Bessel function $J_\nu$. It is  well-known that for $\nu>-1$ the Bessel function $J_\nu(z)$ is positive on the interval  $(0,j_{\nu,1})$. Moreover, $\nu\mapsto j_{\nu,1}$ is increasing, $j_{\nu,1}>\nu$, and $j_{1,1}=3.8317\ldots$, see~\cite[\S10.21]{NIST:DLMF}.

\section{The UIG-model: undirected graphs without loops}
\label{modelA}

We start out by describing the first graph evolution model we want to analyze. At each discrete time step we create a vertex or an edge and each created item gets a time stamp. Precisely, either
(i) a new vertex is created, or
(ii) a new undirected edge is added between two \emph{distinct} existing vertices. The two vertices may already be adjacent. So, loops are forbidden while multiple edges are allowed. We investigate the bivariate generating function that corresponding to all histories generated by this process and keeping track of the number of steps and the number of vertices. 

\bdf
A \emph{history} in the UIG-model is a multigraph obtained by the construction described above, using steps (i) and (ii), such that each vertex and each edge receives a time stamp. 
\edf 
Thus, a history records not only the resulting multigraph but also the order in which its vertices and edges were created. In particular, two histories may have the same underlying multigraph while differing in their time labels.

Note that the time stamps are not labels in the classical sense of analytic combinatorics, as there are restrictions on where to put them. For instance, two vertices are needed to create edges. So, no edge can have time stamp one or two. So, technically, the class of histories is an unlabelled combinatorial class in the sense of~\cite{flajolet2009analytic}, which means that we have to set up ordinary generating functions. To this end, let the variable $t$ mark time (number of steps) and $u$ mark the number of vertices, \ie write 
\[
T(t,u)=\sum_{n\ge0}\sum_{m\ge 0} T_{n,m} t^nu^m
\]
where the coefficient $T_{n,m}=[t^nu^m]T(t,u)$ is the number of histories that can be generated in $n$ steps and have $m$ vertices (and hence $n-m$ edges). A step creating a vertex corresponds to a factor $ut$ in the generating function, as we add a vertex and a step. When creating an edge, we must add a step and then choose two distinct vertices. This amounts to two pointing-and-removing operations and then reverting the removals. Moreover, as the edge is undirected, the order of the vertices must not be taken into account. After all, this yields the following differential equation for the generating function $T(t,u)$: 
\begin{equation*}
T(t,u)=1+u t\,T(t,u)+\frac{u^2 t}{2} \partial_u^2 T(t,u).
\end{equation*}
Equivalently, as a formal power series in $t$,
\begin{equation}\label{product-formula}
T(t,u)=\sum_{m\ge 0}^{\infty}\,\prod_{k=1}^{m}\frac{t\,u}{1-\binom{k}{2}t}.
\end{equation}
This can be seen, as each vertex creation yields a factor $tu$ and after having created $k$ vertices so far, there are $\binom k2$ possibilities for each edge that is created before the creation of the $(k+1)$-th vertex. 

We are primarily interested in asymptotics of
$[t^n]T(t,1)$, as $n\to\infty$. The first few terms can be seen from the expansion 
\[
T(t,1)=1+t+t^2+2t^3+6t^4+25t^5+135t^6+909t^7+\Ord{t^8}.
\] 
\begin{remark}
The $m$-th summand in \eqref{product-formula} has poles at $t=2/(k(k-1))$ for $k=2,\dots,m$. Thus, the poles of $T(t,1)$ accumulate at $t=0$ and hence the radius of convergence of the series $T(t,1)$ is zero unless there are some miraculous cancellations. (Note that for any $z>0$ the sum has positive and negative terms!) But a simple lower bound shows that this is indeed the case: Take $1\le k\le n$ and consider the histories that first create $k$ vertices and then $n-k$ edges. They are easy to count and yield  $[t^n]T(t,1)\ge\binom k2^{n-k}$. Choosing $k=\lfloor\sqrt n\rfloor$ gives 
$\bigl([t^n]T(t,1)\bigr)^{1/n}\to\infty$. To circumvent this problem, the asymptotic analysis is performed at the coefficient level rather than by classical singularity analysis at a dominant singularity. 
\end{remark}
\begin{prop}\label{prop:explicit-coefficients}
For \(n\ge 1\), the coefficients of \(T(t,u)\) are given by
\begin{equation}\label{T-coefficients}
[t^n]T(t,u)
=
\sum_{r=1}^n u^r\sum_{k=1}^r
\frac{
\binom{k}{2}^{\,n-1}
(-1)^{k-r}2^{r-1}(2k-1)
}{
(k+r-1)!(r-k)!
}
\end{equation}
with the convention $0^0=1$. In particular, for \(u=1\), there is a constant $c>0$ such that 
\begin{equation}\label{truncated-bessel-sum}
[t^n]T(t,1)
=
\sum_{k=1}^{2n}
\binom{k}{2}^{\,n-1}
\frac{2k-1}{\sqrt{2}}
J_{2k-1}(2\sqrt{2})
+
\Ord{e^{-cn\log n}}
\end{equation}
where $J_\alpha(z)$ denotes the Bessel function of the first kind. 
\end{prop}
\bpf
We decompose the product into partial fractions. One obtains
\begin{align*}
T(t,u)
&=
1+\sum_{r=0}^\infty u^r\sum_{k=1}^r
\frac{t}{1-\binom{k}{2}t}
\prod_{s=1}^{k-1}
\frac{1}{\binom{k}{2}-\binom{s}{2}}
\prod_{s=k+1}^{r}
\frac{1}{\binom{k}{2}-\binom{s}{2}}.
\end{align*}
Using the Gamma function, this can be rewritten as
\begin{align*}
T(t,u)
&=
1+\sum_{r=0}^\infty u^r\sum_{k=1}^r
\frac{
t(-1)^{k-r}2^{r-1}(2k-1)
}{
\left(1-\binom{k}{2}t\right)
\Gamma(k+r)\Gamma(r+1-k)
}.
\end{align*}
Observe that the $r$-th summand in the original product \eqref{product-formula} is of order $\Ord{t^r}$. Therefore, any term with $r>n$ does not contribute to the coefficient of $t^n$. Extracting the coefficient of $t^n$ from the terms $1\le r\le n$ yields \eqref{T-coefficients}. Setting \(u=1\) and rewriting the resulting expression in terms of
Bessel functions gives
\[
[t^n]T(t,1)
=
\sum_{k\ge1}
\binom{k}{2}^{n-1}
\frac{2k-1}{\sqrt{2}}
J_{2k-1}(2\sqrt{2}).
\]
It remains to show that this sum may be truncated at \(k=2n\) with a
superexponentially small error. For fixed $z$ and $\nu\to\infty$, we use the standard large-order estimate
\[
J_\nu(z)
=
\frac{1}{\Gamma(\nu+1)}
\left(\frac{z}{2}\right)^\nu
\left(1+\Ord{\nu^{-1}}\right).
\]
With $z=2\sqrt{2}$ and $\nu=2k-1$, Stirling's formula gives
\[
J_{2k-1}(2\sqrt{2})
=
\frac{1}{\sqrt{4\pi k}}
\exp\left(-2k\log k+\Ord{k}\right).
\]
On the other hand,
\[
\binom{k}{2}^{n-1}
=
\exp\left(2n\log k+\Ord{n}\right).
\]
Hence the logarithm of the absolute value of the $k$-th summand in the
Bessel representation is bounded by
\[
2n\log k-2k\log k+\Ord{k+n}.
\]
For $k>2n$ this is at most
\[
-k\log k+\Ord{k},
\]
and therefore, for all sufficiently large $n$,
\[
\sum_{k>2n}
\binom{k}{2}^{n-1}
\frac{2k-1}{\sqrt{2}}
\left|J_{2k-1}(2\sqrt{2})\right|
=
\Ord{e^{-cn\log n}}
\]
for some constant $c>0$. This proves \eqref{truncated-bessel-sum}.
\epf
\subsection{Asymptotic analysis}
By \eqref{truncated-bessel-sum}, it remains to estimate
\[
S_n
=
\sum_{k=2}^{2n}
\binom{k}{2}^{n-1}
\frac{2k-1}{\sqrt{2}}
J_{2k-1}(2\sqrt{2})
=
\sum_{k=2}^{2n}\exp(\Phi_n(k)),
\]
where
\begin{align}
\Phi_n(x)
&=
(n-1)\bigl(\log x+\log(x-1)-\log 2\bigr)
+\log(2x-1)-\frac12\log 2
+\log J_{2x-1}(2\sqrt{2}).
\label{phase-definition}
\end{align}
\begin{remark} 
Note that for $n\ge 2$ the first term of the sum in \eqref{T-coefficients} is zero. Thus, we define $S_n$ as sum starting at $k=2$ and avoid the singularity of the logarithm in \eqref{phase-definition}.
\end{remark}
\begin{remark}
Recall the notations and facts on Bessel function presented at the end of Section~\ref{intro}. As $x\ge1$, we have $j_{2x-1,1}\ge j_{1,1}>2\sqrt2$ and hence $J_{2x-1}(2\sqrt2)>0$ for all real $x>1$. Thus $\Phi_n$ is well defined and smooth on $(1,\infty)$ (in fact, the Bessel function $J_\nu(z)$ is analytic in $\nu$). In particular, all summands of $S_n$ are positive.
\end{remark}

We apply the discrete Laplace method to determine the asymptotic behaviour
of \(S_n\) and show the following result.
\begin{theorem}
\label{undirected-asymptotic}
Let $W$ denote the principal Lambert $W$-function
\cite{corless1996lambertw}, defined for $x>0$ by
\[
W(x)e^{W(x)}=x.
\]The number of undirected histories that can be generated in
exactly $n$ steps satisfies
\[
[t^n]T(t,1)
\sim
\sqrt{2}\,
\frac{\LambertW\lp n\sqrt{2}\rp}
     {n2^{2n}\sqrt{1+\LambertW\lp n\sqrt{2}\rp}}
\exp\lp
2n\lp
\LambertW\lp n\sqrt{2}\rp
-1
+\frac{1}{\LambertW\lp n\sqrt{2}\rp}
\rp
\rp.
\]
\end{theorem}
\begin{proof}
Throughout the proof we write
\[
\omega_n:=\LambertW\lp n\sqrt2\rp.
\]
We also set
\[
f(x):=J_{2x-1}(2\sqrt2).
\]

\medskip
\noindent
\emph{Step 1: Estimates for the Bessel factor.}
From the series expansion of the Bessel function we have, for $x$ tending to infinity,
\begin{equation}\label{bessel-series-expansion}
f(x)
=
\sum_{\ell\ge0}
\frac{(-1)^\ell 2^{\ell+x-\frac12}}
     {\ell!\,\Gamma(\ell+2x)}
=
\frac{2^{x-\frac12}}{\Gamma(2x)}
\left(1-\frac1x+\Ord{\frac1{x^2}}\right).
\end{equation}
Consequently,
\[
\log f(x)
=
\left(x-\frac12\right)\log2
-\log\Gamma(2x)
+\Ord{\frac1x}.
\]
Furthermore, by differentiating the series termwise and using the standard
digamma expansion
\[
\frac{\Gamma'(z)}{\Gamma(z)}
=
\log z-\frac{1}{2z}+\Ord{\frac1{z^2}},
\qquad z\to\infty,
\]
we get
\begin{equation}\label{bessel-log-derivative}
\frac{f'(x)}{f(x)}
=
\log2-2\log(2x)+\frac{1}{2x}
+\Ord{\frac{\log x}{x^2}},\qquad \text{as }x\to\infty .
\end{equation}
Differentiating the preceding asymptotic expansion gives, as $x\to\infty$,
\begin{equation}\label{bessel-second-log-derivative}
\left(\frac{f'}{f}\right)'(x)
=
-\frac2x-\frac{1}{2x^2}
+\Ord{\frac{\log x}{x^3}},
\end{equation}
and
\begin{equation}\label{bessel-third-log-derivative}
\left(\frac{f'}{f}\right)''(x)
=
\Ord{\frac1{x^2}}.
\end{equation}

\medskip
\noindent
\emph{Step 2: Location of the saddle point.}
Differentiating \eqref{phase-definition} gives
\[
\Phi_n'(x)
=
(n-1)\left(\frac1x+\frac1{x-1}\right)
+\frac{2}{2x-1}
+\frac{f'(x)}{f(x)}.
\]
Using \eqref{bessel-log-derivative} and expanding in powers of $1/x$,
we obtain
\begin{equation}\label{phase-first-derivative}
\Phi_n'(x)
=
\frac{2n}{x}
-2\log x
-\log2
-\frac{1}{2x}
+\Ord{\frac{n}{x^2}}.
\end{equation}

\emph{Existence and uniqueness of the critical point.}
By \eqref{bessel-second-log-derivative} there is $X_0>1$, independent of $n$, such that $(f'/f)'(x)<0$ for $x\ge X_0$. Since
\[
\Phi_n''(x)
=
-(n-1)\left(\frac1{x^2}+\frac1{(x-1)^2}\right)
-\frac{4}{(2x-1)^2}
+
\left(\frac{f'}{f}\right)'(x),
\]
all terms are negative on $[X_0,\infty)$, so $\Phi_n'$ is strictly decreasing there. 
Note that $f$ is analytic in $\mathbb C$ and positive on the compact interval $[1,X_0]$. So, $f'/f$ is continuous there and hence bounded by some constant $C$. Therefore, 
$\Phi_n'(x)\ge 2(n-1)/X_0-C>0$ for $x\in(1,X_0]$ and all sufficiently large $n$. Finally, $\Phi_n'(x)\to-\infty$ as $x\to\infty$ by \eqref{phase-first-derivative}. Consequently, for all sufficiently large $n$, $\Phi_n'$ has exactly one zero $x^\star\in(X_0,\infty)$; $\Phi_n$ is increasing on $(1,x^\star]$ and decreasing on $[x^\star,\infty)$, and attains its unique maximum at $x^\star$.

\emph{Asymptotic location.}
From the saddle-point equation $\Phi_n'(x)=0$, we obtain the
asymptotic saddle-point equation
\[
\frac{n}{x}
\sim
\log\lp x\sqrt{2}\rp,
\qquad n\to\infty.
\]
To solve this equation, set $y=n/x$. Then
\begin{equation*}
y=\log\left(\frac{n\sqrt{2}}{y}\right),
\qquad\text{and hence}\qquad
y e^y=n\sqrt{2}.
\end{equation*}
Therefore,
\begin{equation*}
y=\LambertW\left(n\sqrt{2}\right)=:\omega_n,
\qquad
x=\frac{n}{\omega_n}.
\end{equation*}
Now let $x^\star$ denote the solution of the saddle-point equation
$\Phi_n'(x)=0$. Writing $v=\log\lp x^\star\sqrt2\rp$,
the saddle-point equation becomes
\[
\frac{n}{x^\star}
=
v+\Ord{\frac1{x^\star}}
+\Ord{\frac{n}{(x^\star)^2}}.
\]
Since $x^\star$ is of order $n/\omega_n$, this implies
\begin{equation*}
\frac{n}{x^\star}=v+\Ord{\frac{\omega_n^2}{n}}.
\end{equation*}
Using \(x^\star=e^v/\sqrt{2}\), we obtain
\[
v e^v
=
n\sqrt2+\Ord{\omega_n}.
\]
By the mean value theorem applied to $z\mapsto ze^z$, and since
\[
\left.\partial_v\left(v e^v\right)\right|_{v=\omega_n}
=
(1+\omega_n)e^{\omega_n}
\asymp n.
\]
we obtain
\[
v
=
\omega_n+\Ord{\frac{\omega_n}{n}}.
\]
Converting back and inserting the expansion for $v$ we get
\[
x^\star  =
\frac{e^{\omega_n}}{\sqrt{2}}
\exp\left(\Ord{\frac{\omega_n}{n}}\right) =
\frac{n}{\omega_n}
\left(1+\Ord{\frac{\omega_n}{n}}\right).
\]
Thus the saddle point satisfies
\[
x^\star
=
\frac{n}{\omega_n}
+\Ord{1}.
\]
We shall also use that this critical point is unique. For
$x\le n^{1/2+\varepsilon}$, with fixed $0<\varepsilon<1/2$, the positive
term $2n/x$ in \eqref{phase-first-derivative} dominates the logarithmic
terms, and therefore $\Phi_n'(x)>0$ for all sufficiently large $n$.
For $x\ge n^{1/2+\varepsilon}$, the right-hand side of
\eqref{phase-first-derivative} is strictly decreasing up to an error which
is smaller than the main derivative. Hence $\Phi_n'$ has exactly one
zero. Thus $\Phi_n$ has a unique maximum at $x^\star$.

\medskip
\noindent
\emph{Step 3: Local quadratic expansion.}
Differentiating once more gives
\[
\Phi_n''(x)
=
-(n-1)\left(\frac1{x^2}+\frac1{(x-1)^2}\right)
-\frac{4}{(2x-1)^2}
+
\left(\frac{f'}{f}\right)'(x).
\]
Using \eqref{bessel-second-log-derivative}, we obtain
\begin{equation}\label{phase-second-derivative}
\Phi_n''(x)
=
-\frac{2n}{x^2}
-\frac2x
+\Ord{\frac{n}{x^3}}
+\Ord{\frac1{x^2}}.
\end{equation}
At the saddle point this yields
\begin{equation}\label{curvature-at-saddle}
\Phi_n''(x^\star)
=
-\frac{2\omega_n(\omega_n+1)}{n}
\left(1+o(1)\right).
\end{equation}
We define
\[
\lambda_n:=-\Phi_n''(x^\star).
\]
Hence
\[
\lambda_n
=
\frac{2\omega_n(\omega_n+1)}{n}
\left(1+o(1)\right).
\]
In particular, the saddle point $x^\star$ corresponds
to a strict local maximum of $\Phi_n$. Similarly,
\[
\Phi_n'''(x)
=
2(n-1)\left(\frac1{x^3}+\frac1{(x-1)^3}\right)
+\frac{16}{(2x-1)^3}
+
\left(\frac{f'}{f}\right)''(x),
\]
and from \eqref{bessel-third-log-derivative} we get, uniformly for
$|x-x^\star|\le\sqrt n$,
\begin{equation}\label{phase-third-derivative}
\Phi_n'''(x)
=
\Ord{\frac{\omega_n^3}{n^2}},\quad \text{as} \ n\to\infty.
\end{equation}
Let
\[
\mathcal W_n
:=
\left\{
k\in\mathbb Z:\ |k-x^\star|\le \sqrt n
\right\}.
\]
Using \eqref{phase-third-derivative}, Taylor's theorem gives, for $x$ with $|x-x^\star|\le\sqrt n$,
\[
\Phi_n(x)
=
\Phi_n(x^\star)
+\frac12\Phi_n''(x^\star)(x-x^\star)^2
+
\Ord{\frac{\omega_n^3}{n^2}|x-x^\star|^3}.
\]
Since $|x-x^\star|\le\sqrt n$, the remainder is
\[
\Ord{\frac{\omega_n^3}{\sqrt n}}=o(1).
\]
Thus, uniformly on $\mathcal W_n$,
\begin{equation}\label{local-quadratic-expansion}
\Phi_n(k)
=
\Phi_n(x^\star)
+\frac12\Phi_n''(x^\star)(k-x^\star)^2
+o(1).
\end{equation}
Exponentiating this expansion yields
\[
e^{\Phi_n(k)}
\sim
e^{\Phi_n(x^\star)}
\exp\left(
-\frac{\lambda_n}{2}(k-x^\star)^2
\right).
\]
The second exponential is referred to as the Gaussian factor arising from the quadratic expansion. It describes the local concentration of the summands around the saddle point.

\medskip\noindent
\emph{Step 4: Negligibility of the tails and the Gaussian sum.}
We split the sum into a central part and a tail part. Then
\[
S_n=S_n^{\mathrm{cent}}+S_n^{\mathrm{tail}} 
\]
where
\[
S_n^{\mathrm{cent}}:= 
\sum_{k\in\mathcal W_n} e^{\Phi_n(k)}   \ \text{ and } \
S_n^{\mathrm{tail}}:=
\sum_{\substack{1\le k\le 2n\\ |k-x^\star|>\sqrt n}}
e^{\Phi_n(k)}.
\]

Since \(\Phi_n\) is increasing on \((1,x^\star)\) and decreasing on
\((x^\star,2n)\), every \(k\) with
\(\lvert k-x^\star\rvert>\sqrt n\) satisfies
\begin{equation*}
\Phi_n(k)
\leq
\max_{\varepsilon\in{-1,1}}
\Phi_n\left(x^\star+\varepsilon\sqrt n\right).
\end{equation*}
Hence, by Taylor's theorem, for some points
\(\xi_\pm\) between \(x^\star\) and \(x^\star\pm\sqrt n\),
\begin{equation*}
\Phi_n\left(x^\star\pm\sqrt n\right)=
\Phi_n(x^\star)
+\frac{n}{2}\Phi_n''(\xi_\pm)
\leq
\Phi_n(x^\star)-c\omega_n^2
\end{equation*}
for some \(c>0\) and all sufficiently large \(n\). Consequently,
\begin{equation*}
S_n^{\mathrm{tail}}
\leq
2n\exp\left(\Phi_n(x^\star)-c\omega_n^2\right)
=e^{\Phi_n(x^\star)}
\cdot\Ord{n e^{-c\omega_n^2}}.
\end{equation*}
Since $\omega_n\sim\log n$, this is negligible compared to the central
contribution. 

It remains to estimate $S_n^{\mathrm{cent}}$. Using \eqref{local-quadratic-expansion}, we get
\[
S_n^{\mathrm{cent}}
=
e^{\Phi_n(x^\star)}
\sum_{k\in\mathcal W_n}
\exp\left(-\frac{\lambda_n}{2}(k-x^\star)^2\right)
(1+o(1)).
\]
The Gaussian tail outside $\mathcal W_n$ is exponentially small, hence the
sum may be extended to all integers:
\[
\sum_{k\in\mathcal W_n}
\exp\left(-\frac{\lambda_n}{2}(k-x^\star)^2\right)
\sim
\sum_{k\in\mathbb Z}
\exp\left(-\frac{\lambda_n}{2}(k-x^\star)^2\right), \text{ as } \nti.
\]
By Poisson summation,
\[
\sum_{k\in\mathbb Z}
\exp\left(-\frac{\lambda_n}{2}(k-x^\star)^2\right)
=
\sqrt{\frac{2\pi}{\lambda_n}}
\left(
1+\Ord{e^{-2\pi^2/\lambda_n}}
\right),
\]
Since $\lambda_n\to0$, this is asymptotically equal to $\sqrt{2\pi / \lambda_n}$.
Consequently,
\[
S_n^{\mathrm{cent}}
=
e^{\Phi_n(x^\star)}
\sqrt{\frac{2\pi}{-\Phi_n''(x^\star)}}
(1+o(1)).
\]
Since \(S_n^{\mathrm{tail}}\) is negligible, we obtain
\begin{equation}\label{laplace-reduction}
S_n
=
e^{\Phi_n(x^\star)}
\sqrt{\frac{2\pi}{-\Phi_n''(x^\star)}}
(1+o(1)).
\end{equation}

\medskip
\noindent
\emph{Step 5: Evaluation of the saddle point and curvature.}
It remains to evaluate the two factors in
\eqref{laplace-reduction}. From Step~1,
\[
\log f(x)
=
\left(x-\frac12\right)\log2
-\log\Gamma(2x)
+\Ord{\frac1x}.
\]
Using Stirling's formula,
\[
\log\Gamma(2x)
=
\left(2x-\frac12\right)\log(2x)
-2x+\frac12\log(2\pi)
+\Ord{\frac1x},
\]
we obtain
\[
\log f(x)
=
-\left(2x-\frac12\right)\log x
+2x
-x\log2
-\frac12\log(2\pi)
+\Ord{\frac1x}.
\]
Substituting this into \eqref{phase-definition} and expanding
\[
\log(x-1)
=
\log x-\frac1x+\Ord{\frac1{x^2}},
\qquad
\log(2x-1)
=
\log x+\log2+\Ord{\frac1x},
\]
gives
\begin{align}
\Phi_n(x)
&=
\left(2n-\frac12\right)\log x
-2x\log x
+2x
-x\log2
-n\log2
+\frac32\log2
-\frac12\log(2\pi)        \notag\\
&\qquad
-\frac{n-1}{x}
+\Ord{\frac{n}{x^2}+\frac1x}.
\label{refined-phase}
\end{align}
By the saddle-point estimate established in Step 2,
\[
x^\star-\frac{n}{\omega_n}=\Ord1.
\]
Since \(\Phi_n'(x^\star)=0\) and
\[
\Phi_n''(x)=\Ord{\frac{\omega_n^2}{n}}
\]
uniformly for $x$ between \(x^\star\) and \(n/\omega_n\), Taylor's theorem gives
\[
\Phi_n(x^\star)=\Phi_n\left(\frac{n}{\omega_n}\right)
+
\Ord{\frac{\omega_n^2}{n}}=
\Phi_n\left(\frac{n}{\omega_n}\right)+o(1).
\]
Evaluating \eqref{refined-phase} at \(x=n/\omega_n\) and using
\[
\log\left(\frac{n}{\omega_n}\right)=\omega_n-\frac12\log2
\]
therefore gives
\begin{equation}\label{phase-at-saddle}
\Phi_n(x^\star)=-2n\log2
+
2n\lp \omega_n-1+\frac1{\omega_n}\rp
-\frac32\omega_n
+\frac54\log2
-\frac12\log\pi
+o(1).
\end{equation}
Using the curvature estimate $\eqref{curvature-at-saddle}$ derived in Step 3, we obtain
\begin{equation}\label{gaussian-prefactor}
\sqrt{\frac{2\pi}{-\Phi_n''(x^\star)}}
=
\sqrt{\frac{\pi n}{\omega_n(\omega_n+1)}}(1+o(1)).
\end{equation}
Combining \eqref{laplace-reduction},
\eqref{phase-at-saddle}, and
\eqref{gaussian-prefactor}, we obtain
\[
S_n
\sim
\exp\left(
-2n\log2
+
2n\left(\omega_n-1+\frac1{\omega_n}\right)
\right)
\exp\left(
-\frac32\omega_n
+\frac54\log2
-\frac12\log\pi
\right)
\sqrt{\frac{\pi n}{\omega_n(\omega_n+1)}}.
\]
Equivalently,
\[
S_n
\sim
2^{-2n}
\exp\left(
2n\left(\omega_n-1+\frac1{\omega_n}\right)
\right)
2^{5/4}
n^{1/2}
e^{-3\omega_n/2}
\frac{1}{\sqrt{\omega_n(\omega_n+1)}}.
\]
Finally, since $\omega_n e^{\omega_n}=n\sqrt2$,
we have
\[
e^{-3\omega_n/2}
=
\left(\frac{\omega_n}{n\sqrt2}\right)^{3/2}.
\]
Therefore
\[
2^{5/4}
n^{1/2}
e^{-3\omega_n/2}
\frac{1}{\sqrt{\omega_n(\omega_n+1)}}
=
\sqrt2\,
\frac{\omega_n}{n\sqrt{1+\omega_n}}.
\]
Hence
\[
S_n
\sim
\sqrt2\,
\frac{\omega_n}{2^{2n}n\sqrt{1+\omega_n}}
\exp\left(
2n\left(\omega_n-1+\frac1{\omega_n}\right) 
\right).
\qedhere
\]
\end{proof}
\section{Distribution of the number of vertices in the UIG-model}
\label{sec:params-undirected}
We now refine the enumeration of histories of the UIG-model by keeping track of the number of
vertices as well. Let $U_n$ denote the number of vertices in a uniformly random
history generated in $n$ steps. Its probability
generating function is
\[
P_n(u)
=
\frac{[t^n]T(t,u)}{[t^n]T(t,1)}.
\]
As in the univariate case, the coefficient admits a similar representation involving Bessel functions, namely, for \(n\ge2\), we have 
\begin{equation}\label{bessel-marked-undirected}
[t^n]T(t,u)
=
\sum_{k\ge2}
\binom{k}{2}^{n-1}
\frac{2k-1}{\sqrt2}
\sqrt{u}\,
J_{2k-1}\bigl(2\sqrt{2u}\bigr),
\end{equation}
for $u$ in a suitable fixed compact neighbourhood $K$ of $1$. 
More precisely, put $u_0:=j_{1,1}^2/8=1.835\ldots$. For $0<u<u_0$ we have $2\sqrt{2u}<j_{1,1}\le j_{2k-1,1}$, so that all summands in \eqref{bessel-marked-undirected} are positive. Throughout this section we fix a compact interval $U$ containing 0 and such that with $e^s<u_0$ for all $s\in U$. We write $u=e^s$. The sum can again be truncated at $2n$, with an error being exponentially smaller than the main term and uniform for $s\in U$.
As in the univariate case we write 
\begin{equation}\label{Laplace-approx-bivariate}
    [t^n]T(t,u)\sim \sum_{k=2}^{2n} \exp(\Phi_{n,u}(k) 
\end{equation}
with 
\begin{equation}\label{def_phi_n_u}
\Phi_{n,u}(x)
:=
\log\left(
\binom{x}{2}^{n-1}
\frac{2x-1}{\sqrt2}
\sqrt{u}\,
J_{2x-1}\bigl(2\sqrt{2u}\bigr)
\right).
\end{equation}
By the series expansion of the Bessel function,
\[
\sqrt u\, J_{2x-1}(2\sqrt{2u}\,)=\frac{2^{x-\frac12}u^x}{\Gamma(2x)}\,\psi(x,s)
\ \text{ with } \ 
\psi(x,s):=\sum_{\ell\ge0}\frac{(-2e^s)^\ell \Gamma(2x)}{\ell!\,\Gamma(\ell+2x)},
\]
we obtain 
\begin{equation}\label{diff-with-remainder}
\Phi_{n,e^s}(x)
=
\Phi_n(x)+sx+R(x,s)
\end{equation}
where 
\[
R(x,s):=\log\psi(x,s)-\log\psi(x,0).
\]
Since $0<\Gamma(2x)/\Gamma(\ell+2x)\le(2x)^{-\ell}$ and
$0\le\partial_x\log\bigl(\Gamma(\ell+2x)/\Gamma(2x)\bigr)\le \ell/x$,
term-wise differentiation shows that $\psi(x,s)=1-e^s/x+\Ord{x^{-2}}$, that
$\partial_s^j\psi(x,s)=-e^s/x+\Ord{x^{-2}}$ for $j\ge1$. Likewise, we obtain 
$\partial_x\partial_s^j\psi(x,s)=\Ord{x^{-2}}$ for $j\ge0$, uniformly for $s\in U$. Consequently,
\begin{equation}\label{R-bounds}
\partial_s^jR(x,s)=\Ord{\frac1x},
\qquad
\partial_x\partial_s^jR(x,s)=\Ord{\frac1{x^2}},
\ \text{ for } \ j=0,1,2,
\end{equation}
as $x\to\infty$, uniformly for $s\in U$.
\begin{prop}
\label{prop:undirected-marked-saddle}
For $s\in U$ and $n$ sufficiently large, the equation
$\Phi'_{n,e^s}(x)=0$ has a unique solution $x_n(s)$ in $(1,\infty)$; it is a smooth function of $s$.
Set 
\(
\omega_n(s):=\LambertW\bigl(n\sqrt2\,e^{-s/2}\bigr).
\) Then, uniformly for $s\in U$,
\[
x_n(s)
=
\frac{n}{\omega_n(s)}
+\Ord1.
\]
\end{prop}
\begin{proof}
Existence and uniqueness follow verbatim as in Step~2 of the proof of Theorem~\ref{undirected-asymptotic}. Indeed, by \eqref{diff-with-remainder} and \eqref{R-bounds} the function $\Phi_{n,e^s}$ satisfies the same estimates as $\Phi_n$, uniformly for $s\in U$; smoothness in $s$ follows from the implicit function theorem, as $\partial_x^2\Phi_{n,e^s}(x_n(s))<0$.
Together with \eqref{phase-first-derivative}, \eqref{diff-with-remainder} and \eqref{R-bounds} yield
\[
\Phi'_{n,e^s}(x)
=
\frac{2n}{x}
-2\log x
-\log2
+s
+\Ord{\frac{n}{x^2}}
+\Ord{\frac1x}.
\]
Hence the saddle-point equation $\Phi'_{n,e^s}(x)=0$ implies
\[
\frac{n}{x_n(s)}
=
\log\!\left(x_n(s)\sqrt2\,e^{-s/2}\right)
+
\Ord{\frac{n}{x_n(s)^2}}
+
\Ord{\frac1{x_n(s)}}.
\]
As in the unmarked case, this shows uniformly for $s$ near $0$ that
\[
x_n(s)\asymp\frac{n}{\log n}, \ \text{ as } \nti.
\]

Set
\[
v:=\log\!\left(x_n(s)\sqrt2\,e^{-s/2}\right).
\]
Multiplying the saddle-point equation by $x_n(s)$ and using
$x_n(s)\asymp n/\log n$, we obtain
\[
v e^v
=
n\sqrt2\,e^{-s/2}
+
\Ord{\omega_n(s)}, \ \text{ as }\nti.
\]
On the other hand,
\[
\omega_n(s)e^{\omega_n(s)}
=
n\sqrt2\,e^{-s/2}.
\]
Since $s$ stays in the fixed neighbourhood $U$ of $0$,
we have the uniform estimate $\omega_n(s)\asymp\log n$. The mean value theorem applied to
$v\mapsto ve^v$ therefore gives
\[
v=\omega_n(s)+\Ord{\frac{\omega_n(s)}{n}}.
\]
Consequently, for $\nti$, we obtain 
\[
x_n(s)
=
\frac{e^{s/2}}{\sqrt2}e^v
=
\frac{n}{\omega_n(s)}
\exp\!\left(\Ord{\frac{\omega_n(s)}{n}}\right)
=
\frac{n}{\omega_n(s)}
+\Ord1.
\qedhere\]
\end{proof}
\begin{theorem}[Gaussian limit law for the number of vertices]
\label{undirected-vertices-normal-law}
Let $U_n$ denote the number of vertices in a uniformly random history of the UIG-model generated in $n$ steps. With
\[
\omega_n:=\LambertW(n\sqrt2),
\]
we have
\[
\EE[U_n]\sim \frac{n}{\omega_n},
\qquad
\Var(U_n)\sim \frac{n}{2\omega_n(\omega_n+1)}.
\]
Moreover,
\[
\frac{U_n-\EE[U_n]}{\sqrt{\Var(U_n)}}
\dlim
\mathcal N(0,1).
\]
\end{theorem}

\begin{proof}
We study the moment generating function
\[
P_n(e^s)
=
\frac{[t^n]T(t,e^s)}{[t^n]T(t,1)}
\]
for $s$ in a fixed neighbourhood of $0$ and write
\[
K_n(s):=\log P_n(e^s).
\] 
By virtue of the uniform bounds on the remainder $R(x,s)$ established in \eqref{R-bounds}, the discrete Laplace argument used in the univariate enumeration applies uniformly for $s$ near $0$. This yields
\[
[t^n]T(t,e^s)
=
\exp\bigl(\Phi_{n,e^s}(x_n(s))\bigr)
\sqrt{\frac{2\pi}{-\Phi_{n,e^s}''(x_n(s))}}
(1+o(1)).
\]
To justify the extraction of moments, let $G_n(s) := \Phi_{n,e^s}(x_n(s))$ denote the dominant saddle-point exponent. We must verify that differentiating the exact sum representation \eqref{bessel-marked-undirected} of $[t^n]T(t,e^s)$ with respect to $s$ is asymptotically equivalent to differentiating its dominant term. Recall from \eqref{Laplace-approx-bivariate} that the $k$-th summand is $\exp(\Phi_{n,e^s}(k))$. Differentiating this with respect to $s$ brings a factor 
\[
\partial_s \Phi_{n,e^s}(k) = k + \partial_s R(k,s) = x_n(s) + (k - x_n(s)) + \partial_s R(k,s).
\]
Repeated differentiation with respect to $s$ thus introduces polynomial factors in $(k-x_n(s))$, as well as derivatives of $R(k,s)$. Thanks to the uniform bounds provided by \eqref{R-bounds}, the derivatives of the remainder are asymptotically negligible. The polynomial factors $(k-x_n(s))^j$ cannot counter the strong Gaussian decay within the central window. As the summands remain exponentially small outside this window, the derivatives of $K_n(s)$ up to third order are asymptotically equivalent to the corresponding derivatives of $G_n(s)$. Thus, we can evaluate them by applying the multivariate chain rule directly to $G_n(s)$, which yields 
\[
G_n'(s) = \partial_x \Phi_{n,e^s}(x_n(s)) x_n'(s) + \partial_s \Phi_{n,e^s}(x_n(s)).
\]

By definition, the saddle point satisfies
\[
\Phi'_{n,e^s}(x_n(s))=0,
\]
Note that despite $\Phi_{n,e^s}$ has $s$ in its index, we can rely on the representation given in \eqref{diff-with-remainder} for the computation. Doing so, we obtain
\[
K_n'(s)=x_n(s)+\Ord1,
\qquad
K_n''(s)=x_n'(s)+\Ord1,
\qquad
K_n'''(s)=\Ord{\frac{n}{\omega_n(s)^3}}.
\]
Here, in the first derivative, the term involving $x_n'(s)$ vanishes, because
$x_n(s)$ is the saddle point.

Likewise, we determine the variance by differentiation of the the saddle-point equation with respect to \(s\) and relying on \eqref{diff-with-remainder}. Hence, we get
\[
\left(\partial_{xx}\Phi_{n,e^s}\right)(x_n(s))x_n'(s)
+
\left(\partial_{xs}\Phi_{n,e^s}\right)(x_n(s))=
0,
\]
and hence
\[
x_n'(s)=-
\frac{
\left(\partial_{xs}\Phi_{n,e^s}\right)(x_n(s))
}{
\left(\partial_{xx}\Phi_{n,e^s}\right)(x_n(s))
}.
\]
From \eqref{diff-with-remainder} we directly obtain 
\[
\left(\partial_{xs}\Phi_{n,e^s}\right)(x_n(s))=1+o(1),
\]
whereas the refined curvature estimate from the univariate saddle-point
analysis yields
\[
-\left(\partial_{xx}\Phi_{n,e^s}\right)(x_n(s))=\frac{
2\omega_n(s)\bigl(\omega_n(s)+1\bigr)
}{n}
\bigl(1+o(1)\bigr).
\]
Consequently,
\[
x_n'(s)=\frac{n}{
2\omega_n(s)\bigl(\omega_n(s)+1\bigr)
}
\bigl(1+o(1)\bigr).
\]
Evaluating at $s=0$ now gives
\[
\EE[U_n]= K_n'(0)=\frac{n}{\omega_n}+\Ord1
\]
and
\[
\Var(U_n)=
K_n''(0)=\frac{n}{2\omega_n(\omega_n+1)}
\bigl(1+o(1)\bigr).
\]
For the limit law, set
\[
\mu_n:=\EE[U_n],
\qquad
\sigma_n^2:=\Var(U_n).
\]
For fixed $\tau\in\mathbb R$, let $h_n=\tau/\sigma_n$. Taylor's theorem
gives
\[
K_n(h_n)=\mu_n h_n
+
\frac{\sigma_n^2h_n^2}{2}
+
\frac{K_n'''(\xi_n)}{6}h_n^3
\]
for some $\xi_n$ between $0$ and $h_n$. Since
\[
K_n'''(\xi_n)=\Ord{\frac{n}{\omega_n^3}}
\qquad\text{and}\qquad
\sigma_n^2
\asymp
\frac{n}{\omega_n^2},
\]
the remainder is $\Ord{n^{-1/2}}=o(1)$. Therefore,
\[
\log
\EE\exp\left(
\tau\frac{U_n-\mu_n}{\sigma_n}
\right)=\frac{\tau^2}{2}+o(1).
\]
The moment generating functions therefore converge in a neighbourhood of
the origin to that of the standard normal distribution, which proves the
claim.
\end{proof}
\begin{remark}
The proof gives $\EE[U_n]=n/\omega_n+o(\sqrt n/\omega_n)$. Numerical computations suggest that in fact $\EE[U_n]=n/\omega_n+\Ord1$.
\end{remark}

\section{A weighted version of the UIG-model}
\label{sec:weighted-undirected}
In this section we are interested in what changes if one introduces weights on the creation of vertices and edges. In particular, phase transitions can be expected here, but a totally rigorous treatment of these would require a delicate saddle-point analysis that is uniform in the weights. This is beyond the scope of this paper and so we confine ourselves with a treatment of the case of constant weights and discuss expected phenomena briefly. 

The parameter analysis in the previous section shows that in the unweighted model the number of vertices is typically of order
\[
\frac{n}{\LambertW(n\sqrt2)}\sim \frac{n}{\log n}.
\]
Thus the proportion of vertex-creation steps tends to zero, and asymptotically almost all steps are edge-creation steps. It is
therefore natural to introduce a relative weight for vertex creation as a way to tune the vertex--edge balance. We assign weight $\alpha > 0$ to a
vertex-creation step and weight $\beta > 0$ to an edge-creation step, obtaining
\[
T_{\alpha,\beta}(t,u)
=
1+\alpha u t\,T_{\alpha,\beta}(t,u)
+
\beta\frac{u^2t}{2}
\partial_u^2T_{\alpha,\beta}(t,u).
\]
We may set $\tau=\beta t$ and replace $\alpha$ by $\alpha/\beta$. Since this only contributes a global factor $\beta^n$
to coefficients of $t^n$, it cancels in all conditional distributions with $n$ steps. Hence, without loss of generality, we take $\beta=1$
from now on and write simply $\alpha$ for the relative vertex weight. Thus, no new coefficient extraction is needed: the weighted model is obtained
from the previous bivariate series by the substitution $u\mapsto \alpha u$.

For $\alpha>0$, let $U_n^{(\alpha)}$ denote the number of vertices
in a random history chosen with probability proportional to $\alpha^m$, where $m$ is its number of vertices. Its probability generating
function is
\[
P_{n,\alpha}(u)
=
\frac{[t^n]T(t,\alpha u)}
     {[t^n]T(t,\alpha)}.
\]
For each fixed $\alpha$, the same saddle-point argument applies with the
marking parameter centred at $s=\log\alpha$. Since $\alpha e^s=e^{\log\alpha+s}$, the weighted model is obtained by shifting the logarithmic marking parameter by $\log\alpha$. Thus the relevant saddle at $s=0$ is the saddle of the original model at $s=\log\alpha$. 

Since $\alpha e^s=e^{\log\alpha+s}$, the weighted model is obtained by shifting the logarithmic marking parameter by $\log\alpha$. Thus the relevant saddle-point at $s=0$ in the weighted model corresponds to a saddle at $s=\log\alpha$ in the original model. 
Therefore, for fixed $\alpha\in(0,u_0)$, the results of Section~\ref{sec:params-undirected} apply verbatim, with $U$ replaced by a compact interval containing $\log\alpha$ and contained in $(-\infty,\log u_0)$. For fixed $\alpha\ge u_0$, finitely many summands of \eqref{bessel-marked-undirected} (those with $2k-1\le 2\sqrt{2\alpha}$, since $j_{\nu,1}>\nu$) may be nonpositive; their total contribution is $\Ord{k_0^{2n}}$ for some fixed $k_0$, which is negligible compared with the main term, and the analysis restricted to $k>k_0$ goes through unchanged. We omit the details.
Hence, with
\[
\omega_{n,\alpha}
:=
\LambertW\left(n\sqrt{\frac{2}{\alpha}}\right),
\]
we obtain
\[
\EE\left[U_n^{(\alpha)}\right]
\sim
\frac{n}{\omega_{n,\alpha}},
\qquad
\Var\left[U_n^{(\alpha)}\right]
\sim
\frac{n}{2\omega_{n,\alpha}(\omega_{n,\alpha}+1)}
\sim
\frac{n}{2\omega_{n,\alpha}^2}.
\]
Moreover, the same logarithmic moment-generating-function argument gives
asymptotic normality after centering and normalization.

As long as $\alpha>0$ stays fixed, we still have
$\omega_{n,\alpha}\sim\log n$, and consequently
\[
\EE\left[U_n^{(\alpha)}\right]
\sim \frac{n}{\log n}.
\]
Thus, for every such $\alpha$ the vertex proportion still tends to zero. The weight changes the lower-order location and variance of the vertex count, but not its leading term $n / \log n$.

Finally, it is instructive to consider heuristically what happens when
$\alpha=\alpha_n$ is allowed to grow with $n$. The saddle-point expression
suggests that, as long as $\alpha_n=o(n^2)$, the dominant number of
vertices remains sublinear. A crossover is expected when $\alpha_n$ is of order $n^2$. The Lambert $W$ term is then of constant order, and the saddle becomes linear in $n$.

This is coherent with the fact that, for a generic weight $\alpha$, the recurrence relation for the coefficients \(
a_{n,k}^{(\alpha)} := [t^n u^k]T(t,\alpha u)\)
becomes
\begin{equation*}
    a^{(\alpha)}_{n+1,k+1}
    =
    \alpha a^{(\alpha)}_{n,k}
    +
    \frac{k(k+1)}{2}a^{(\alpha)}_{n,k+1}.
\end{equation*}
The vertex-creation weight is of order $\alpha$, whereas the number of possible edge-creation steps is quadratic in the number of vertices. Thus, when the relevant number of vertices is of order $n$, the scale $\alpha_n\asymp n^2$ is the natural one at which the balance between vertex and edge creation may change. A rigorous treatment of this regime would require a saddle-point analysis uniform in $\alpha_n$. But this at the first glance innocent looking modification seems to be a Herculean task, as we expect that this requires taking into account the more error terms in all the computations the closer 
$\alpha_n$ is to the order $n^2$.

\section{The DIG-model: directed graphs with loops allowed}
\label{sec:directed}
We now consider the directed variant of the previous model. At each
discrete time step, either
(i) a new vertex is created, or
(ii) a new directed edge is added between two existing vertices. In this
model loops are allowed. Thus, if there are currently \(k\) vertices,
there are \(k^2\) possible directed edges that may be created. We again refer to the resulting time-labelled multigraph as a history.

Let \(t\) mark time and let \(u\) mark the number of vertices. A vertex
creation contributes a factor \(ut\). Edge creation corresponds to choosing
an ordered pair of vertices, possibly the same vertex twice, and is represented by the operator
$
u  \partial_u\!\left(ut\,\partial_u\right)
$.
Therefore the generating function $\overrightarrow{T}(t,u)$ satisfies
\begin{equation}\label{directed-functional-equation}
\overrightarrow{T}(t,u)
=
1
+
ut\,\overrightarrow{T}(t,u)
+
u\partial_u\!\left(
ut\,\partial_u\overrightarrow{T}(t,u)
\right).
\end{equation}
The coefficient \([t^n u^m]\overrightarrow{T}(t,u)\) is the number of
directed graph histories generated in \(n\) steps with \(m\)
vertices. We are interested in the asymptotic behaviour of
\[
[t^n]\overrightarrow{T}(t,1),
\qquad \text{as } n\to\infty.
\]
At \(u=1\), the first few terms are
\[
\overrightarrow{T}(t,1)
=
1+t+2t^2+7t^3+37t^4+264t^5+2433t^6+27913t^7+\Ord{t^8}.
\]
We start with an auxiliary lemma. 
\begin{lemma}\label{auxlem} 
For integers $m>\ell\ge 1$ the following identity holds: 
\[
\sum_{k=1}^m \frac{2 (-1)^{m-k}}{(m-k)!(m+k)!} k^{2\ell} = 0. 
\]
\end{lemma}
\bpf
Consider the function $f(z)=\prod_{j=1}^m 1/(z-j^2)$ and write it as a partial fraction decomposition $f(z)=\sum_{k=1}^m C_k/(z-k^2)$. Then, computing the coefficients with the help of residues gives 
\[
C_k = \prod_{1\le j\le m, \  j \neq k} \frac{1}{k^2 - j^2} = \frac{2 k^2 (-1)^{m-k}}{(m-k)!(m+k)!}.
\]
Now expand $1/(z-k^2)$ into a geometric series for large $z$. We get 
\[
f(z)=\sum_{k=1}^m \sum_{r\ge 0} \frac{C_k k^{2r}}{z^{r+1}}
= \sum_{r\ge 0}\sum_{k=1}^m\frac{C_k k^{2r}}{z^{r+1}}
\]
and thus
\[
[z^{-\ell}] F(z)= \sum_{k=1}^m C_k k^{2\ell-2} = \sum_{k=1}^m \frac{2 (-1)^{m-k}}{(m-k)!(m+k)!} k^{2\ell}. 
\]
As $f(z)=\Ord{z^-m}$ for $z\to\infty$, the inner sum must be zero for $\ell<m$, as claimed. 
\epf

Now we turn back to the generating function of our graph histories. The next result shows that the radius of convergence is zero, but a certain Borel-like transform of it has an explicit and simple shape. 

\begin{prop}
The bivariate generating function $\overrightarrow{T}(t,u)$ of the number of histories in the DIG-model keeping track of the number of vertices as well has radius of convergence zero, when regarded as power series in $t$ and can be expressed as  
\begin{equation}\label{directed-product-formula}
\overrightarrow{T}(t,u)
=
\sum_{m=0}^{\infty}
\prod_{k=1}^{m}
\frac{tu}{1-k^2t}.
\end{equation}
Let 
\begin{equation}\label{directed-borel-coefficient-relation}
\mathcal{B}\left(\sum_{n\ge0}a_n t^n\right)
:=
\sum_{n\ge0}a_n\frac{(2t)^n}{(2n)!}
\end{equation}
and define 
$\overrightarrow{B}(t,u):=\mathcal{B}\bigl(\overrightarrow{T}(\cdot,u)\bigr)(t)$, \ie
we apply $\mathcal B$ to $\overrightarrow{T}(t,u)$ regarded as power series in $t$. 
Then we have the explicit closed form
\begin{equation}\label{directed-borel-closed-form}
\overrightarrow{B}(t,u)
=
\cE\Big(u\big(\cE(t)-1\big)\Big),
\ \text{ where } \ 
\cE(y):=\cosh \big(\sqrt{2y}\big).
\end{equation}
\end{prop}
\begin{remark}
In general, the Borel transform of order $\kappa$ is defined by
\[
\widehat{\mathcal B}_\kappa
\left(\sum_{n\ge0}a_n t^n\right)
:=
\sum_{n\ge0}
\frac{a_n}{\Gamma(1+n/\kappa)}t^n,
\]
see~\cite{Balser1994,MichalikTkacz2019}. For $\kappa=1/2$,
$\Gamma(1+n/\kappa)=\Gamma(1+2n)=(2n)!$. 
Thus, $\mathcal B$ is the Borel transform of order $1/2$ followed by the rescaling $t\mapsto 2t$.
\end{remark}

\begin{remark}
The formula \eqref{directed-product-formula} suggests that the radius of convergence is zero. As in the UIG-model, there may occur cancellations in the sum such not all poles of a term are actual poles of $\overrightarrow{T}(t,u)$. But a simple lower bound for the number of histories can be obtained by considering the histories that first create $k$ vertices and then $n-k$ edges. This gives $[t^n]\overrightarrow{T}(t,1)\ge k^{2(n-k)}$. The choice $k=\lfloor\sqrt n\rfloor$ shows that the radius of convergence is indeed zero. 
\end{remark}

\begin{proof}
We get \eqref{directed-product-formula} directly from the $k^2$ possible edge creation steps between the creation of the $k$-th and ($k+1$)-st vertices.
A partial fraction decomposition gives
\[
\prod_{k=1}^m \rcp{1-k^2t} = 2\sum_{k=1}^m (-1)^{m-k} \frac{k^{2m}}{(m-k)!(m+k)!}\cdot\frac{1}{1-k^2t}, 
\]
since the residue coefficients are $\prod_{j\ne k}k^2/(k^2-j^2)=k^{2m-2}C_k$ with $C_k$ as in the proof of Lemma~\ref{auxlem}.
Thus 
\begin{align*}
\overrightarrow{T}(t,u)=&1+ 2\sum_{m\ge1}\sum_{k=1}^m \frac{(-1)^{m-k}u^m}{(m-k)!(m+k)!}
\sum_{n\ge 0} k^{2(n+m} t^{n+m} \\
=&1+2\sum_{\ell\ge 1} t^\ell \sum_{m=1}^\ell u^m \sum_{k=1}^m k^{2\ell} \frac{(-1)^{m-k}}{(m-k)!(m+k)!} \\
=&1+2\sum_{\ell\ge 1} t^\ell \sum_{m\ge1} u^m \sum_{k\ge1} k^{2\ell} \frac{(-1)^{m-k}}{(m-k)!(m+k)!} \\
=&1+2\sum_{\ell\ge 1} t^\ell \sum_{k\ge1}^\infty k^{2\ell} \sum_{m\ge k}^\infty  \frac{(-1)^{m-k} u^{(m-k)+k}}{(m-k)!((m-k)+2k)!} \\
=&1+2\sum_{\ell\ge 1} t^\ell \sum_{k\ge1}^\infty k^{2\ell} J_{2k}(2\sqrt u\,)
\end{align*}
where the extension of the summation range of $m$ to $m>\ell$ in the third line is justified by Lemma~\ref{auxlem}, and the interchange of summations in the fourth line by absolute convergence (for fixed $\ell$ and $u$).

Next we apply the rescaled Borel transform $\mathcal B$ and obtain 
\begin{align*}
\overrightarrow{B}(t,u)=& 1+2\sum_{\ell\ge 1} \frac{2t)^\ell}{(2\ell)!}\sum_{k\ge1}^\infty k^{2\ell} J_{2k}(2\sqrt u\,) \\
=&1 + 2\sum_{k\ge } J_{2k}(2\sqrt{u}\,) \sum_{\ell \ge1}\frac{(k\sqrt{2t})^{2\ell}}{(2\ell)!} \\
=& 1 + 2\sum_{k\ge } J_{2k}(2\sqrt{u}\,) \left( \cosh(k\sqrt{2t}\,) - 1 \right) \\
=& 1 + J_0(2\sqrt{u}\,) + 2\sum_{k\ge1} J_{2k}(2\sqrt{u}\,) \cosh(k\sqrt{2t}\,) - \left( J_0(2\sqrt{u}\,) + 2\sum_{k\ge1} J_{2k}(2\sqrt{u}\,) \right).
\end{align*}

Now, we use the generating function
\[
\sum_{n\in\mathbb{Z}}J_n(z)v^n
=
\exp\left(\frac{z}{2}\left(v-\frac{1}{v}\right)\right),
\qquad v\neq0,
\]
see~\cite[Eq.~10.12.1]{NIST:DLMF}, together with $J_{-n}=(-1)^nJ_n$. Putting $v=e^\theta$ and taking the part that is even in $\theta$ gives
\[
J_0(z)+2\sum_{k\ge1}J_{2k}(z)\cosh(2k\theta)=\cosh(z\sinh\theta);
\]
in particular, for $\theta=0$, $J_0(z)+2\sum_{k\ge1}J_{2k}(z)=1$. With $z=2\sqrt u$ and $\theta=\sqrt{2t}/2$ we obtain
\[
\overrightarrow{B}(t,u)=\cosh\left(2\sqrt u\,\sinh\frac{\sqrt{2t}}{2}\right)=\cosh\left(\sqrt{2u\bigl(\cosh\sqrt{2t}-1\bigr)}\right),
\]
because $\cosh\sqrt{2t}-1=2\sinh^2(\sqrt{2t}/2)$. This is \eqref{directed-borel-closed-form}, as desired.
\end{proof}

\begin{remark}
This closed form is structurally analogous to the exponential generating
function \(\exp(u(e^t-1))\) for set partitions: both occurrences of the
exponential function are replaced by
\(\cE(y)=\cosh(\sqrt{2y})\).
More precisely, the coefficients
\([t^n]\overrightarrow{T}(t,1)\) form the even-index subsequence
of the central Bell numbers associated with the central factorial
numbers of the second kind~\cite{KimKimCentralBell}. The same sequence
is listed in OEIS A135920, where it is referred to as a generalized
Bell number sequence. A combinatorial explanation of
this analogy, in terms of partitions into blocks of even size, is given
in Section~\ref{subsec:directed-even-partitions}.
\end{remark}
\begin{remark}
The (exponential generating) function in~\eqref{directed-borel-closed-form} has also been studied in the literature, but in a different context. It is obviously differentially algebraic, but its ordinary counterpart, (its inverse Borel transform) is not, as Klazar~\cite{klazar2003} has shown.    
\end{remark}
To obtain coefficient asymptotics for
$\overrightarrow{B}(t,1)$, we use Hayman's theory of
$\mathcal H$-admissible functions~\cite{Hayman1956}. They are amenable to the saddle-point method in a straight-forward fashion by design. For such functions, say $f(z)$, 
the saddle point $\rho_n$ is determined via the logarithmic derivative
$a_f(r):=r f'(r)/f(r)$ 
by $a_f(\rho_n)=n$, while the second-order logarithmic derivative $b_f(r):=ra_f'(r)$
governs the width of the Gaussian-like peak around the saddle point. Once $\mathcal H$-admissibility is established, these quantities yield the standard Hayman coefficient
asymptotics~\cite{Hayman1956} based on the saddle-point method. We therefore first collect the corresponding data needed for
$\overrightarrow{B}(t,1)$, which will later also be used in the bivariate
analysis for $u$ near $1$.
\begin{lemma}[Hayman data]\label{hayman-data-directed}
For simplicity, set
\(
g(t):=\cE(t)-1
\), with $E(t)$ given in \eqref{directed-borel-closed-form}, and define for any function $h(z)$ its first-order and second-order logarithmic derivatives by  
\begin{equation} \label{log_deriv}
a_h(z):=\frac{z\,h'(z)}{h(z)},
\qquad
b_h(z):=za_h'(z).
\end{equation} 
Then
\[
a_{\cE}(z)
=
\frac{\sqrt{2z}}{2}\tanh(\sqrt{2z}),
\qquad
b_{\cE}(z)
=
\frac{\sqrt{2z}}{4}\tanh(\sqrt{2z})
+
\frac{z}{2}\operatorname{sech}^{2}(\sqrt{2z}),
\]
and
\[
a_g(z)
=
\sqrt{\frac z2}\coth\left(\sqrt{\frac z2}\right),
\qquad
b_g(z)
=
\frac12\sqrt{\frac z2}\coth\left(\sqrt{\frac z2}\right)
-
\frac{z}{4}\operatorname{csch}^{2}\left(\sqrt{\frac z2}\right).
\]
Moreover, as \(z\to\infty\), we have 
\[
\cE(z)\sim g(z)\sim \frac12 e^{\sqrt{2z}},
\qquad
a_{\cE}(z)\sim a_g(z)\sim \sqrt{\frac z2},
\qquad
b_\cE(z)\sim b_g(z)\sim \sqrt\frac z8.
\]
\end{lemma}
\begin{proof}
By \eqref{directed-borel-closed-form} we have
\[
\cE'(z)
=
\frac{\sinh(\sqrt{2z})}{\sqrt{2z}},
\]
and therefore
\[
a_{\cE}(z)
=
\frac{\sqrt{2z}}{2}\tanh(\sqrt{2z}).
\]
Differentiating once more and using
$b_{\cE}(z)=za_{\cE}'(z)$ gives
\[
b_{\cE}(z) = \frac{\sqrt{2z}}{4}\tanh(\sqrt{2z})
+ \frac{z}{2}\operatorname{sech}^2(\sqrt{2z}).
\]
The standard hyperbolic identities then give 
\[
a_g(z)
=
\sqrt{\frac z2}\coth\left(\sqrt{\frac z2}\right) \]
and therefore
\[
b_g(z) = za_g'(z) = \frac12\sqrt{\frac z2}\coth\left(\sqrt{\frac z2}\right)
- \frac{z}{4}\operatorname{csch}^2\left(\sqrt{\frac z2}\right).
\]
The stated asymptotic relations follow immediately from
$\tanh x\to 1$, $\coth x\to1$, $\sech x\to0$, $\csch x\to0$ and $\cosh x\sim \frac12e^x$, as $x\to\infty$. 
\end{proof}
\begin{theorem}[Asymptotic number of histories]
\label{directed-asymptotic}
The entire function \(t\mapsto \overrightarrow{B}(t,1)\) given by
\eqref{directed-borel-closed-form} is $\mathcal H$-admissible in the sense of~\cite{Hayman1956} in the whole complex plane. Let $w_n=\LambertW(2n)$, then, as \(n\to\infty\),
\[
c_n;=[z^n]\overrightarrow{T}(t,1)
\sim
\frac{1}{\sqrt{w_n+1}}\;
\frac{
\exp\!\left(2n\left(w_n+\frac{1}{w_n}\right)\right)
}
{(2e)^{2n}}.
\]
\end{theorem}

\begin{proof}
With $E(t)$ and $g(t)$ as in Lemma~\ref{hayman-data-directed}, write
\[
F(t):=\overrightarrow{B}(t,1)=\cE(g(t)).
\]
Since \(\cE\) and \(g\) are entire and \(g(r)>0\) for \(r>0\),
the function \(F\) is entire and \(F(r)>0\) for \(r>0\).
We verify Hayman admissibility in the sense of
\cite{Hayman1956}. To this end, first write $a(z):=a_F(z)$ and $b(z):=b_F(z)$. We have to verify three conditions: First, the second-order logarithmic derivative $b(z)$ must tend to infinity on the positive real axis. Second, $F(z)$ has an asymptotic representation as $|z|\to\infty$ involving the Hayman data and which is valid near the real axis. Precisely, there is a function $\delta$ such that for $z=re^{i\theta}$ with $r\to \infty$ the asymptotic relation 
\[
F(re^{i\theta}) \sim F(r) e^{ia(r)\theta-b(r)\theta^2/2}
\]
holds uniformly for $|\theta|\le \delta(r)$, \ie the modulus of $F(z)$ is Gaussian near the real axis. Third, farther away from the real axis the modulus of $F(z)$ is small and there is an upper bound involving only the Hayman data. Precisly, we must have $F(z) = o (F(r) / \sqrt{b(r)}\,)$, as $r\to \infty$, uniformly for $|\theta|\ge \delta(r)$. By Lemma~\ref{hayman-data-directed} and the chain rule, the Hayman data are
\begin{equation}\label{univariate-hayman-data}
a(r)
=
a_g(r)\,a_{\cE}(g(r)),
\qquad
b(r)
=
b_g(r)\,a_{\cE}(g(r))
+
a_g(r)^2\,b_{\cE}(g(r)).
\end{equation}
For defining the region around the real axis, fix
\[
\delta(r):=b(r)^{-2/5}.
\]

\medskip
\noindent
\textbf{Condition~\cite[(2.6)]{Hayman1956}: divergence of \(b(r)\).}
Using \eqref{univariate-hayman-data} and Lemma~\ref{hayman-data-directed}, we have

$$
a_g(r)\sim \sqrt{\frac r2},
\qquad
b_g(r)\sim \sqrt{\frac r8},
\qquad
g(r)\sim \frac12 e^{\sqrt{2r}}.
$$
Moreover,
$$
a_{\cE}(z)\sim \sqrt{\frac z2},
\qquad
b_{\cE}(z)\sim \sqrt{\frac z8},
$$
and hence
$$
a_{\cE}(g(r))
=
\Theta\left(e^{\sqrt{r/2}}\right),
\qquad
b_{\cE}(g(r))
=
\Theta\left(e^{\sqrt{r/2}}\right).
$$
Therefore the two terms in \eqref{univariate-hayman-data} satisfy
$$
b_g(r)a_{\cE}(g(r))
=
\Theta\left(\sqrt r\,e^{\sqrt{r/2}}\right),
$$
whereas
$$
a_g(r)^2b_{\cE}(g(r))
=
\Theta\left(r\,e^{\sqrt{r/2}}\right).
$$
Thus the second term dominates, and consequently
$$
b(r)
\sim
a_g(r)^2\,b_{\cE}(g(r))
=
\Theta\left(r\,e^{\sqrt{r/2}}\right)
\xrightarrow[r\to\infty]{}\infty.
$$
In particular, \(\delta(r)\to0\).

\medskip
\noindent
\textbf{Condition~\cite[(2.4)]{Hayman1956}: local expansion for
\(|\theta|\leq\delta(r)\).}
Set
\[
L_r(\theta):=\log F(r e^{i\theta}).
\]
Writing \(y=\sqrt{r/2}\), we have
\[
F(re^{i\theta})
=
\cosh\left(
2\sinh\left(ye^{i\theta/2}\right)
\right).
\]
Since
$y\delta(r)\longrightarrow0$,
uniformly for \(|\theta|\le\delta(r)\),
\[
\Re\left(
2\sinh\left(ye^{i\theta/2}\right)
\right)
=
2\sinh\left(y\cos(\theta/2)\right)
\cos\left(y\sin(\theta/2)\right)
>0
\]
for all sufficiently large \(r\). The zeros of \(\cosh z\) are purely
imaginary, so \(F(re^{i\theta})\neq0\) throughout this arc. Hence
\(L_r(\theta)=\log F(re^{i\theta})\) can be chosen analytic on a
neighbourhood of the arc.

On this zero-free neighbourhood, let
\[
a(z):=z\frac{F'(z)}{F(z)},
\qquad
b(z):=z a'(z),
\]
so that these functions extend the quantities \(a(r)\) and \(b(r)\)
defined above. Taylor's theorem with remainder yields, uniformly for
\(|\theta|\le \delta(r)\),
\[
L_r(\theta)
=
L_r(0)
+
i\theta a(r)
-
\frac{\theta^2}{2}b(r)
+
\Ord{
|\theta|^3
\sup_{|\varphi|\le|\theta|}
\big|r e^{i\varphi} b'(r e^{i\varphi})\big|
}.
\]
Using the explicit expressions from Lemma~\ref{hayman-data-directed},
differentiation of $b(z)$ from \eqref{univariate-hayman-data} gives, uniformly for
$|\varphi|\le\delta(r)$,
\[
\bigl|
re^{i\varphi}b'(re^{i\varphi})
\bigr|
=
\Ord{\sqrt r\,b(r)}.
\]
Indeed, since \(y\delta(r)=o(1)\), the standard hyperbolic estimates used
in Lemma~\ref{hayman-data-directed} hold uniformly on this shrinking arc,
and hence
\[
a_g(z),\,b_g(z),\,zb_g'(z)=\Ord{\sqrt r},
\]
whereas
\[
a_{\cE}(g(z)),\,b_{\cE}(g(z)),
\,g(z)b_{\cE}'(g(z))
=
\Ord{e^{\sqrt{r/2}}}.
\]
Since $b(r)=\Theta(r e^{\sqrt{r/2}})$, the Taylor remainder is therefore
uniformly bounded by
\[
\Ord{\delta(r)^3\sqrt r\,b(r)}
=
\Ord{\sqrt r\,b(r)^{-1/5}}
=
o(1).
\]
This proves the local Gaussian approximation.

\medskip
\noindent
\textbf{Condition~\cite[(2.5)]{Hayman1956}: tail bound for
\(\delta(r)\le|\theta|\le\pi\).}
Let \(z=r e^{i\theta}\), with \(\delta(r)\le|\theta|\le\pi\), and set
\begin{equation}\label{eta}
\eta(z):=\sqrt{2\lp\cosh\lp\sqrt{2z}\rp-1\rp}.
\end{equation}
Since \(\cosh\) is even, the value of \(\cosh(\eta(z))\) is independent
of the branch of the square root, and $F(z)=\cosh(\eta(z))$.

\emph{Small angles: \(\delta(r)\le|\theta|\le1\).}
Here
\[
\Re \lp\sqrt{2z}\rp=\sqrt{2r}\cos(\theta/2)\to\infty,
\]
and therefore
\[
\cosh \lp\sqrt{2z}\rp-1
=
\frac12 e^{\sqrt{2z}}
\left(1+\Ord{e^{-\Re\lp\sqrt{2z}\rp}}\right).
\]
Consequently,
\[
|\eta(z)|
\sim
\exp\left(\sqrt{r/2}\cos(\theta/2)\right).
\]
Using \(|\cosh w|\le e^{|w|}\), we get
\[
\log|F(z)|
\le
|\eta(z)|
\sim
\exp\left(\sqrt{r/2}\cos(\theta/2)\right).
\]
Since
\[
\log F(r)\sim \exp\left(\sqrt{r/2}\right)
\]
and $1-\cos(\theta/2)\ge \theta^2/12$ 
for $|\theta|\le1$, we get the estimate
\begin{align*}
\log|F(r e^{i\theta})|-\log F(r)&
\le e^{\sqrt{r/2}} \lp \exp\lp -\frac{\theta^2}{12}\sqrt\frac r2\rp-1\rp \\
&\le e^{\sqrt{r/2}} \lp \exp\lp -\frac{\delta(r)^2}{12}\sqrt\frac r2\rp-1\rp
\end{align*}
Now observe $\delta(r)^2\sqrt{r}\asymp r^{-3/10}\exp\lp -(4/5)\sqrt{r/2}\rp\to 0$, as $r\to\infty$, there exists a constant \(c>0\) such that
\begin{equation}\label{minor-arc-argument}
\log|F(r e^{i\theta})|-\log F(r) \le -c r^{-3/10} \exp\lp \frac{\sqrt{r}}{5\sqrt2}\rp =o(\log b(r)) 
\end{equation}
where we used $b(r)\asymp r e^{\sqrt{r/2}}$. 
Hence
\[
|F(r e^{i\theta})|
=
o\left(\frac{F(r)}{\sqrt{b(r)}}\right)
\ 
\text{ for }\ \delta(r)\le|\theta|\le1.
\]

\emph{Large angles: \(1\le|\theta|\le\pi\).}
Then
\[
\Re \lp\sqrt{2z}\rp
=
\sqrt{2r}\cos(\theta/2)
\le
\sqrt{2r}\cos(1/2),
\]
and therefore $|\eta(z)|=\Ord{\exp(\sqrt{r/2}\cos(1/2))}$, which implies
\[
\log|F(z)|
=
\Ord{e^{\sqrt{r/2}\cos(1/2)}}.
\]
Since \(\log F(r)\asymp e^{\sqrt{r/2}}\), this gives $\log|F(z)|=o(\log F(r))$ and thus 
\begin{equation}\label{large-angles}
\log\left(\frac{|F(z)|\sqrt{b(r)}}{F(r)}\right)
\le
\frac12\log b(r)
-
\log F(r)+o(\log F(r))
\xrightarrow[r\to\infty]{}
-\infty.
\end{equation}
This proves the tail condition. Hence \(F(t)=\overrightarrow{B}(t,1)\)
is \(\mathcal H\)-admissible.

What remains to be done in order to derive the asymptotic value of $c_n$ is to determine the asymptotic location of the saddle point $\rho_n$ and of $b(\rho_n$. Writing \(y=\sqrt{r/2}\), we have
\[
a(r)
=
y\cosh(y)\tanh(2\sinh y).
\]
The right-hand side is strictly increasing from \(0\) to \(+\infty\) for
\(y>0\). Hence, for every \(n\ge1\), there is a unique \(\rho_n>0\)
satisfying \(a(\rho_n)=n\).

By Hayman's coefficient theorem~\cite[Corollary~II, p.~69]{Hayman1956},
letting \(\rho_n\) be this unique solution, we have
\[
[t^n]F(t)
\sim
\frac{F(\rho_n)}
     {\rho_n^n\sqrt{2\pi b(\rho_n)}}.
\]
Using the coefficient relation \eqref{directed-borel-coefficient-relation}, this gives
\begin{equation}\label{directed-hayman-transfer}
c_n
\sim
\frac{(2n)!}{2^n}
\frac{F(\rho_n)}
     {\rho_n^n\sqrt{2\pi b(\rho_n)}}.
\end{equation}
It remains to evaluate the terms in \eqref{directed-hayman-transfer}.
Write
\[
y_n:=\sqrt{\frac{\rho_n}{2}}.
\]
Then
\[
g(\rho_n)
=
\cosh(2y_n)-1
=
2\sinh^2 y_n,
\]
and hence
\[
a_g(\rho_n)
=
y_n\coth y_n,
\qquad
a_{\cE}(g(\rho_n))
=
\sinh y_n\,\tanh(2\sinh y_n).
\]
Therefore the saddle-point equation \(a(\rho_n)=n\) becomes
\[
n
=
y_n\cosh(y_n)\tanh(2\sinh y_n).
\]
Using
\[
\cosh y_n
=
\frac{e^{y_n}}{2}\left(1+e^{-2y_n}\right)
\]
and
\[
\tanh(2\sinh y_n)
=
1+\Ord{e^{-4\sinh y_n}},
\]
the saddle-point equation yields
\[
2n
=
y_n e^{y_n}
\left(1+\Ord{e^{-2y_n}}\right).
\]
Comparing the saddle-point equation with the defining relation for $w_n$, \(w_n e^{w_n}=2n,\) and
using
\[
\partial_y\left(ye^y\right)=(1+y)e^y.
\]
the mean value theorem gives directly
\[
y_n
=
w_n+\Ord{e^{-2w_n}}.
\]
Consequently,
\[
\rho_n
=
2w_n^2
\left(
1+\Ord{\frac{w_n}{n^2}}
\right),
\]
and in particular
\[
\rho_n^n
=
(2w_n^2)^n
\exp\left(
\Ord{\frac{w_n}{n}}
\right).
\]
For the evaluation of \(F(\rho_n)\), we use
\[
F(\rho_n)
=
\cosh(2\sinh y_n).
\]
The saddle-point equation implies
\[
2\sinh y_n
=
\frac{2n}{w_n}+o(1),
\]
and therefore
\begin{equation}\label{directed-F-at-saddle}
F(\rho_n)
=
\frac12
\exp\left(\frac{2n}{w_n}\right)
(1+o(1)).
\end{equation}
Next,  using the Hayman data in \eqref{univariate-hayman-data} at \(r=\rho_n=2y_n^2\), we obtain
\[
b(\rho_n)
=
\frac{n}{2}(y_n+1)(1+o(1))
=
\frac{n}{2}(w_n+1)(1+o(1)).
\]
Thus
\begin{equation}\label{directed-b-at-saddle}
\sqrt{2\pi b(\rho_n)}
=
\sqrt{\pi n(w_n+1)}(1+o(1)).
\end{equation}
Finally, Stirling's formula gives
\[
(2n)!
=
\sqrt{4\pi n}\left(\frac{2n}{e}\right)^{2n}(1+o(1)).
\]
Together with the preceding estimate for $\rho_n^n$ and
$w_n e^{w_n}=2n$, we find
\[
\frac{(2n)!}{2^n\rho_n^n}
=
\sqrt{4\pi n}\,
\frac{\exp(2nw_n)}{2^{2n}e^{2n}}
(1+o(1)).
\]
Combining this estimate with \eqref{directed-F-at-saddle} and
\eqref{directed-b-at-saddle} in \eqref{directed-hayman-transfer}, we obtain
\[
c_n
\sim
\frac{1}{\sqrt{w_n+1}}\,
\frac{
\exp\left(2n\left(w_n+\frac1{w_n}\right)\right)
}
{(2e)^{2n}}.
\]
This proves the theorem.
\end{proof}

\begin{remark}
Since $(2e)^{-2n}e^{2nw_n}=2^{-2n}e^{2n(w_n-1)}$, the exponential growth in Theorems~\ref{undirected-asymptotic} and~\ref{directed-asymptotic} has the same shape $2^{-2n}\exp\bigl(2n(\omega-1+1/\omega)\bigr)$, with $\omega=\LambertW(n\sqrt2)$ and $\omega=\LambertW(2n)$, respectively.
\end{remark}

\section{The number of vertices in the DIG-model}
\label{param_B}

In this section we turn to the distribution of the number of vertices in a random history in the DIG-model that is generated in $n$ steps. We will work with the (rescaled) Borel transform $\overrightarrow{B}(t,u)$ of the bivariate generating function and show that like in the UIG-model a central limit theorem holds in the DIG-model. Recall that we set $c:n=[t^n]\overrightarrow{T}(t,1)$ in Theorem~\ref{directed-asymptotic}, Moreover, let $c_{n,m}=[t^nu^m]\overrightarrow{T}(t,u)$ be the number of histories that are generated in $n$ steps and yield a graph with $m$ vertices. Our proof is based on the concept of e-admissibility, see \cite{DrmotaGittenbergerKlausner2005}, which extends Hayman admissibility to functions in two variables and yields a normal limit law. For the sake of the paper being self-contained, we state the definition of e-admissibility. This is essentially~\cite[Def.~1]{DrmotaGittenbergerKlausner2005}, but specialized to entire functions which are positive on the positive real axis. 

\bdf\label{def_e-adm}
Let $f(z,u)$ be a function and let
and let $a(z,u)$, $b(z,u)$ be its first-order and second-order logarithmic derivative with respect to $z$, respectively, and $\bar a(z,u)$, 
$\bar b(z,u)$ their analogues with respect to $u$ and $c(z,u)=u (\partial_u a)(z,u)$ its mixed second-order logarithmic derivative. Finally, let $B(z,u)$ denote the logarithmic Hessian of $f(z,u)$. 

The function $f(z,u)$ is called e-admissible of the following ten conditions are satisfied:
\begin{enumerate}
\item $f(z,u)$ is entire and $f(r,1)>0$ for positive $r$
\item For each $K>0$, setting $\eps(r)= K\sqrt{b(r,1)/\det B(r,1)}=K(\bar b(r,1)-c(r,1)^2/b(r,1))^{-1/2}$, there is a function 
$\delta:(0,\infty)\to(0,\pi)$ such that, as $r\to\infty$, we have 
\[
f(re^{i\theta},u)\sim f(r,u)\exp\left(i\theta
a(r,u)-\frac{\theta^2}2b(r,u)\right),
\]
uniformly for $|\theta|\le\delta(r)$ and $u\in[1-\eps(r),1+\eps(r)]$.
\item 
The estimate  
\[
f(re^{i\theta},u)= o\left(\frac{f(r,u)}{\sqrt{b(r,u)}}\right),
\quad\mbox{as }r\to \infty,
\]
holds uniformly for $\delta(r)\le|\theta|\le\pi$ and $u\in[1-\eps(r),1+\eps(r)]$.
\item \label{bbquer}
For $r\to \infty$ we have $b(r,1)\to\infty$.
\item \label{b-ex}
$b(r,u)\sim b(r,1)$ for $r\to\infty$, uniformly for 
$u\in[1-\eps(r),1+\eps(r)]$. 
\item \label{a-ex}
$a(r,u)=a(r,1)+c(r,1)(u-1)+\Ord{c(r,1)(u-1)^2}$ for $r>0$ and
$u\in[1-\eps(r),1+\eps(r)]$.
\item \label{ad_lemma4}
$\bar a(r,u)=\Ord{\bar a(r,1)}$ and $\bar b(r,u)=\Ord{\bar b(r,1)}$
for all $u$ in an arbitrary but fixed complex neighbourhood of $1$ and all $r>0$.
\item \label{det}
$(\det B(r,1))/b(r,1)\to \infty$ as $r\to \infty$. 
\item\label{eps}
$\eps(r)^3\bar b(r,1)\to 0$ for $r\to\infty$.
\item\label{f}
$\bar b(r,1)=\Ord{\bar a(r,1)^2}$ and $\bar a(r,1)=\Ord{f(r,1)^\eta}$ for every
$\eta>0$, as $r\to\infty$. 
\end{enumerate}
\edf

As shown in~\cite{DrmotaGittenbergerKlausner2005}, a sequence of random variables that is defined via the coefficients of an e-admissible function in the way as in Theorem~\ref{directed-vertices-normal-law} below is asymptotically normally distributed. Like in the UIG-model, the asymptotic expectation and variance can be expressed by means of the Hayman data of the function.  

Now we are in the position to state and prove the central limit theorem for the number of vertices in a random history of given length in the directed model. We will prove it by showing that $\overrightarrow{B}(t,u)$ is e-admissible. 

\begin{theorem}[Number of vertices in the directed model]
\label{directed-vertices-normal-law}
Let \(M_n\) be the number of vertices in a uniformly random history generated in \(n\) steps, \ie
\[
\mathbb P(M_n=m)=\frac{c_{n,m}}{c_n}
\qquad(m\ge0).
\]
Let \(w_n=\LambertW(2n).\)
Then
\[
\mathbb E[M_n]\sim \frac{n}{w_n},
\qquad
\Var(M_n)\sim \frac{n}{2w_n(w_n+1)}.
\]
Moreover,
\[
\frac{M_n-\mathbb E[M_n]}{\sqrt{\Var(M_n)}}
\dlim
\mathcal N(0,1).
\]
\end{theorem}
\begin{proof}
We verify the e-admissibility conditions required by~\cite[Theorem~2]{DrmotaGittenbergerKlausner2005} for
$\overrightarrow{B}(t,u)$ and then apply this theorem to establish
the distributional asymptotics. As noted by the authors, for ~\cite[Theorem~2]{DrmotaGittenbergerKlausner2005} it suffices
to check conditions (1)--(6), (8), and (9) of Definition~\ref{def_e-adm}.

By \eqref{directed-borel-coefficient-relation}, the probability
generating function of $M_n$ is
\[
P_n(u)
=
\frac{[t^n]\overrightarrow{B}(t,u)}
     {[t^n]\overrightarrow{B}(t,1)}.
\]
Thus the distribution of $M_n$ can be studied directly through the
coefficients of $\overrightarrow{B}(t,u)$. In view of Definition~\ref{def_e-adm}, we extend the Hayman quantities $a(r)$ and $b(r)$ introduced in
\eqref{univariate-hayman-data} by defining
\[
a(r,u):=r\partial_r\log\overrightarrow{B}(r,u),
\qquad
b(r,u):=r\partial_r a(r,u).
\]
and additionally set
\[
\bar a(r,u):=u\partial_u\log\overrightarrow{B}(r,u),
\qquad
\bar b(r,u):=u\partial_u\bar a(r,u),
\qquad
c(r,u):=u\partial_u a(r,u).
\]
Using \eqref{directed-borel-closed-form} and the chain rule, we obtain,
with $q=ug(r)$,
\[
\begin{aligned}
a(r,u)
&=a_g(r)a_{\cE}(q),\qquad
b(r,u)=b_g(r)a_{\cE}(q)+a_g(r)^2b_{\cE}(q),\\
\bar a(r,u)&=a_{\cE}(q), \quad\bar b(r,u)
=b_{\cE}(q),\quad c(r,u)=a_g(r)b_{\cE}(q).
\end{aligned}
\]
The variance quantity appearing in Definition~\ref{def_e-adm} is
\[
V(r):=\frac{\det B(r,1)}{b(r,1)}=  
\bar b(r,1)-\frac{c(r,1)^2}{b(r,1)}.
\]
Substituting the preceding identities gives
\[
V(r)=
\frac{
b_g(r)a_{\cE}(g(r))b_{\cE}(g(r))
}{
b_g(r)a_{\cE}(g(r))
+a_g(r)^2b_{\cE}(g(r))
}.
\]
Put $y=\sqrt{r/2}$. Using the asymptotics established in the proof of
Theorem~\ref{directed-asymptotic}, the preceding identity gives
\[
V(r)
\sim
\frac{e^y}{4(y+1)}
\longrightarrow\infty,
\]
which establishes condition (8).

Condition (1) follows from \eqref{directed-borel-closed-form}:
$\overrightarrow{B}(t,u)$ is entire in both variables, its
coefficients are nonnegative, and
$\overrightarrow{B}(r,1)>0$ for $r>0$.

Fix $K>0$ and set
\[
\varepsilon(r):=K V(r)^{-1/2}.
\]
Hence,
\[
\varepsilon(r)
=
\Ord{\sqrt y\,e^{-y/2}}
\longrightarrow0.
\]
Thus, for sufficiently large $r$, the interval
$[1-\varepsilon(r),1+\varepsilon(r)]$ is contained in
$U:=[1/2,3/2]$.

Condition (4) was already established in the proof of
Theorem~\ref{directed-asymptotic}, since $b(r,1)=b(r)$.

For condition (5), the asymptotics in
Lemma~\ref{hayman-data-directed} give, uniformly for $u\in U$,
\[
b(r,u)
\sim
\frac{\sqrt u}{4}y(y+1)e^y.
\]
Hence
\[
\frac{b(r,u)}{b(r,1)}
=
\sqrt u\,(1+o(1))
\]
uniformly for $u\in U$. Since $\varepsilon(r)\to0$, we have
$\sqrt u=1+o(1)$ uniformly for
$|u-1|\le\varepsilon(r)$. Therefore
\[
b(r,u)\sim b(r,1)
\]
uniformly on this interval, which proves condition~(5).

Furthermore, the explicit expression for $a(r,u)$ gives
\[
\sup_{u\in U}|\partial_{uu}a(r,u)|
=
\Ord{ye^y},
\qquad
c(r,1)\asymp ye^y.
\]
Taylor's theorem therefore yields
\[
a(r,u)
=
a(r,1)+c(r,1)(u-1)
+
\Ord{c(r,1)(u-1)^2},
\]
uniformly for $|u-1|\le\varepsilon(r)$, proving condition (6).

Finally, since
$\bar b(r,1)\asymp e^y$,
\[
\varepsilon(r)^3\bar b(r,1)
=
\Ord{y^{3/2}e^{-y/2}}
=o(1),
\]
which establishes condition (9).

It remains to verify the parameter-uniform versions of the local
Gaussian approximation and the tail bound established in
Theorem~\ref{directed-asymptotic}.

For condition~(2), we use the same cutoff
\[
\delta(r)=b(r,1)^{-2/5}
\]
as in Theorem~\ref{directed-asymptotic}. The only additional point is
uniformity in $u$. Writing
\[
H_\theta
=
2\sqrt u\,\sinh\left(ye^{i\theta/2}\right),
\]
we have, uniformly for $u\in U$ and $|\theta|\le\delta(r)$,
\[
\Re H_\theta\asymp e^y,
\qquad
\partial_\theta^3
\log\overrightarrow B(re^{i\theta},u)
=
\Ord{y^3e^y}.
\]
Since
\[
\delta(r)^3y^3e^y
=
\Ord{y^{3/5}e^{-y/5}}
=o(1),
\]
the Taylor argument from Theorem~\ref{directed-asymptotic} applies
uniformly in $u$, proving condition (2).

For the tail bound, condition (3), the argument leading to \eqref{minor-arc-argument} in the proof of Theorem~\ref{directed-asymptotic} is uniform for $u\in U$: indeed, the only difference is a factor $\sqrt u$ in \eqref{eta}, which is bounded from above and away from zero in $U$. Thus the computations leading to \eqref{minor-arc-argument} and \eqref{large-angles} apply almost verbatim and yield the desired tail bound. 
We have therefore verified all conditions required by
Theorem~2 of~\cite{DrmotaGittenbergerKlausner2005}. This implies so far that $M_n$ converges to a normal limiting distribution. 

\medskip
\noindent
What remains to do is to determine asymptotic values for expectation and variance. Let $\rho_n$ and $y_n=\sqrt{\rho_n/2}$ be as in the proof of
Theorem~\ref{directed-asymptotic}. There we showed that
\[
y_n=w_n+\Ord{e^{-2w_n}}.
\]
Using the saddle-point equation, we obtain
\[
\bar a(\rho_n,1)
=
\sinh y_n\,\tanh(2\sinh y_n)
=
\frac{n\tanh y_n}{y_n}
\sim\frac{n}{w_n}.
\]
Moreover, 
\[
V(\rho_n)
\sim
\frac{n}{2y_n(y_n+1)}
\sim
\frac{n}{2w_n(w_n+1)}.
\]
By Theorem~2 of
\cite{DrmotaGittenbergerKlausner2005},
\[
\EE[M_n]
=
\bar a(\rho_n,1)+o(V(\rho_n)),
\qquad
\Var(M_n)\sim V(\rho_n).
\]
Hence, we finally obtain 
\[
\EE[M_n]\sim\frac{n}{w_n},
\qquad
\Var(M_n)\sim
\frac{n}{2w_n(w_n+1)},
\]
and
\[
\frac{M_n-\EE[M_n]}{\sqrt{\Var(M_n)}}
\dlim
\mathcal N(0,1).
\]
as claimed in the assertion. 
\end{proof}

\begin{remark}
In order to show asymptotic normality we relied on~\cite[Theorem~2]{DrmotaGittenbergerKlausner2005}, for which only eight of the ten properties defining e-admissibility are needed. The other two are only important to prove the closure properties of the class of e-admissible functions. Therefore we confined ourselves with proving only those eight properties for our function. But it can be shown that $\overrightarrow{B}(t,u)$ satisfies the remaining two conditions as well and is therefore indeed e-admissible. 
\end{remark}

\section{A correspondence between the DIG-model and set partitions}
\label{subsec:directed-even-partitions}

As already mentioned in this paper, the function $\overrightarrow{B}(t,u)$ has a striking similarity to the generating function of set partitions. We simply have to replace the standard exponential function $e^t$ by $\cE(t)$. Thus, the coefficients of $\overrightarrow{T}(t,u)$ are closely related to partitions into blocks of even size.

\begin{prop}
For $n \geq 1, 1 \leq k \leq n$ let $p_{n,k}$ denote the number of
partitions of the $2n$-element set $\{1,1',2,2',\ldots,n,n'\}$ into $k$ blocks, each of even cardinality. 
\begin{equation}\label{directed-even-partition-correspondence}
[t^n u^k]\overrightarrow{T}(t,u)
=
\frac{p_{n,k}}{(2k-1)!!}
\end{equation}
where $(2k-1)!!=(2k-1)(2k-3)\cdots3\cdot1$.
\end{prop}

\begin{remark}
    This relation is recorded in OEIS A156289. For the sake of being self-contained, we present a proof here. 
\end{remark}
\begin{proof}
The exponential generating function of the numbers $p_{n,k}$ is
\[
\sum_{n\ge0} p_{n,k}\frac{z^{2n}}{(2n)!}
=
\frac{(\cosh z-1)^k}{k!}.
\]
Recall that $g(t)=\cE(t)-1=\cosh(\sqrt{2t})-1$. If we set
$z=\sqrt{2t}$, this becomes
\[
\sum_{n\ge0} p_{n,k}\frac{(2t)^n}{(2n)!}
=
\frac{g(t)^k}{k!}
\]
with the conventions
\[
p_{0,0}=1,\qquad p_{n,0}=0\quad(n\ge1).
\]
On the other hand, the closed form~\eqref{directed-borel-closed-form} for the rescaled Borel transform gives
\[
[u^k]\overrightarrow{B}(t,u)
=
\frac{2^k}{(2k)!}g(t)^k
=
\frac{1}{(2k-1)!!}\frac{g(t)^k}{k!}.
\]
Comparing coefficients with the definition of the rescaled Borel transform
yields \eqref{directed-even-partition-correspondence}.
\end{proof}

The numbers $p_{n,k}$ satisfy the recurrence relation
\cite{comtet2012advanced}
\begin{equation}\label{even-partition-recurrence}
p_{n,k}
=
(2k-1)p_{n-1,k-1}
+
k^2p_{n-1,k}.
\end{equation}
This also follows directly from the generating function. Indeed, if
\[
P_k(z):=\frac{(\cosh z-1)^k}{k!},
\]
then
\[
P_k''(z)
=
(2k-1)P_{k-1}(z)+k^2P_k(z),
\]
and extracting the coefficients gives
\eqref{even-partition-recurrence}.

This recurrence admits a combinatorial interpretation based on the
insertion of the new pair $\{n,n'\}$ into a partition of a $(2n-2)$-element set. Throughout, we assume that the ground set is totally ordered by
\[
1\prec1'\prec2\prec2'\prec\cdots\prec n\prec n'.
\]
and that the blocks of a partition receive a time stamp when they are created. 

Let $\pi$ be a partition of $\{1,1',\ldots,n-1,(n-1)'\}$ into even-sized blocks. The new pair is inserted by one of the following operations: 
\begin{itemize}
\item[(N0)] $\{n,n'\}$ forms a new block.
\item[(N1)] Choose an existing block $B$ from $\pi$ and $\varepsilon\in\{n,n'\}$, put $\varepsilon$ into $B$ and the other element $\eta$ of $\{n,n'\}$ into a new block. Then move $x:=\max_\prec B$ (the maximum of the \emph{old} elements of $B$) to the new block. The result is the replacement of $B$ by the two blocks $\{\eta,x\}$ and  $(B\setminus\{x\})\cup\{\varepsilon\}$ in $\pi$. 
\item[(L)] Choose an existing block $B$ and put both new elements $n$ and $n'$ into $B$.
\item[(E)] Choose an ordered pair $(B,B')$ of distinct existing blocks and put $n$ into $B$ and $n'$ into $B'$. Then move $x:=\max_\prec(B\cup B')$ (taken over the old elements) from the block containing it to the other one.
\end{itemize}
If $\pi$ has $k-1$ blocks, then there are $1+2(k-1)=2k-1$ possibilities to apply either of the operations (N0) and (N1), each resulting in a partition with $k$ blocks. 
If $\pi$ has $k$ blocks, then there are $k+k(k-1)=k^2$ possibilities to apply either (L) or (E). All operations produce partitions into even-sized blocks. 

After normalization by $(2k-1)!!$, recurrence
\eqref{even-partition-recurrence} becomes
\[
\frac{p_{n,k}}{(2k-1)!!}
=
\frac{p_{n-1,k-1}}{(2k-3)!!}
+
k^2\frac{p_{n-1,k}}{(2k-1)!!}.
\]
Using \eqref{directed-even-partition-correspondence}, this is precisely
\[
[t^n u^k]\overrightarrow{T}(t,u)
=
[t^{n-1}u^{k-1}]\overrightarrow{T}(t,u)
+
k^2[t^{n-1}u^k]\overrightarrow{T}(t,u),
\]
which exactly matches the recurrence for the DIG-model on the coefficient level. Thus, the factor
$(2k-1)!!$ is the normalization that transforms the recurrence for
even-block partitions into the recurrence for directed histories. 
In both cases, the two terms correspond, respectively, to the creation of a new vertex and to the $k^2$ possible directed edge or loop choices.

The operations of the insertion description above are reversible. From the resulting partition, one can uniquely recover the preceding partition and the operation that was applied. Iterating the construction therefore yields a correspondence with histories in the DIG-model, after all. At the $j$-th vertex-creation step, there are exactly $2j-1$ possible partition insertions corresponding to the same vertex-creation step, whereas an edge or loop step determines its insertion uniquely.

Consequently, every history in the DIG-model with $k$ vertices has exactly
\[
\prod_{j=1}^{k}(2j-1)=(2k-1)!!
\]
preimages among the even-sized set partitions.

\section{Uniform sampling of graphs from the DIG-model}
\label{algorithm}

In order to sample histories in the DIG-model such that all histories generated in $n$ steps are equally likely, we use the framework of Boltzmann samplers~\cite{Duchon2004}. However, as the generating function $\overrightarrow{T}(z,1)$ is only a formal power series, this does not work directly. Likewise, the Borel transform $\overrightarrow{B}(t,1)$ cannot be used directly in the classical way of the symbolic method, as its building block $E(t)$ does not correspond to any combinatorial construction presented in~\cite{flajolet2009analytic}. 

A naive way of sampling histories in the DIG-model would be by uniformly generating set partitions into blocks of even size and then applying the above correspondence to obtain a graph. However, this procedure does not produce a uniform distribution over fixed length histories of the DIG-model. The reason is that every graph with $k$ vertices originates from $(2k-1)!!$ partitions and $k$ is not constant over all histories generating in $n$ steps. To correct this bias, we present a biased sampler for set partitions into even-sized blocks with the probability of each partition being biased by a factor $1/(2k-1)!!$. When finally applying the correspondence, all histories generated in $n$ steps are equally likely.

\begin{figure}
    \centering
    \resizebox{\textwidth}{0.85\textheight}{
    \begin{tikzpicture}[
    >=Stealth,
    partition/.style={rectangle, draw=black, rounded corners, fill=gray!10, thick, minimum width=8cm, minimum height=0.6cm, align=center, font=\small},
    temp/.style={rectangle, draw=red, dashed, fill=red!5, thick, minimum width=8cm, minimum height=0.6cm, align=center, font=\small},
    history/.style={rectangle, draw=black, fill=yellow!15, thick, align=center, font=\footnotesize\ttfamily, inner sep=4pt},
    action/.style={font=\footnotesize\itshape, text=black!80, align=left},
    graphnode/.style={circle, draw=black, thick, fill=green!10, minimum size=0.7cm, inner sep=1pt, font=\small},
    graphedge/.style={->, thick, draw=blue!80},
    title/.style={font=\large\bfseries}
]

\node[font=\bfseries] at (0, 0.7) {Partition State ($P_m$)};
\node[font=\bfseries] at (7.9, 0.7) {History Stack};

\node[partition] (P5) at (0, 0) {$B_1=\{1, 2'\} \quad B_2=\{1', 2, 4, \mathbf{5}\} \quad B_3=\{3, 3', 4', \mathbf{5'}\}$};

\draw[->, thick] (P5) -- (0, -1.5) node[midway, right, action] {Remove $\mathbf{5, 5'}$ (Diff blocks)};
\node[temp] (T5) at (0, -1.5) {$B_1=\{1, 2'\} \quad B_2=\{1', 2, 4\} \quad B_3=\{3, 3', \mathbf{4'}\}$};

\draw[->, thick] (T5) -- (0, -3.2) node[midway, right, action] {Parity Restore: Transfer max ($4'$) from $B_3$ to $B_2$ \\ $B_3$ not empty \fra{(case (E))}};
\node[partition] (P4) at (0, -3.2) {$B_1=\{1, 2'\} \quad B_2=\{1', 2, \mathbf{4}, \mathbf{4'}\} \quad B_3=\{3, 3'\}$};

\node[history] (H5) at (7.9, -1.5) {\begin{tabular}{c} \textbf{EDGE} $v_2 \to v_3$ \end{tabular}};
\draw[->, dashed, text=gray] (4.0, -1.5) -- (H5) node[midway, above, font=\tiny, text=gray] {Push};

\draw[->, thick] (P4) -- (0, -5) node[midway, right, action] {Remove $\mathbf{4, 4'}$ (Same block: \fra{case (L)}) \\ No parity restore needed};
\node[partition] (P3) at (0, -5) {$B_1=\{1, 2'\} \quad B_2=\{1', 2\} \quad B_3=\{\mathbf{3}, \mathbf{3'}\}$};

\node[history] (H4) at (7.9, -4) {\begin{tabular}{c} \textbf{LOOP} at $v_2$ \\ \hline \textbf{EDGE} $v_2 \to v_3$ \end{tabular}};
\draw[->, dashed, text=gray] (5.6, -4) -- (H4) node[midway, above, font=\tiny, text=gray] {Push};

\draw[->, thick] (P3) -- (0, -7) node[midway, right, action] {Remove $\mathbf{3, 3'}$ \\ $B_3$ becomes empty (\fra{case (N1)})};
\node[partition] (P2) at (0, -7) {$B_1=\{1, \mathbf{2'}\} \quad B_2=\{1', \mathbf{2}\}$};

\node[history] (H3) at (7.9, -6) {\begin{tabular}{c} \textbf{NODE} $v_3$ \\ \hline \textbf{LOOP} at $v_2$ \\ \hline \textbf{EDGE} $v_2 \to v_3$ \end{tabular}};
\draw[->, dashed, text=gray] (2.2, -5.95) -- (H3) node[midway, above, font=\tiny, text=gray] {Push};

\draw[->, thick] (P2) -- (0, -8.5) node[midway, right, action] {Remove $\mathbf{2, 2'}$ (Diff blocks)};
\node[temp] (T2) at (0, -8.5) {$B_1=\{1\} \quad B_2=\{\mathbf{1'}\}$};

\draw[-, thick] (T2) -- (0, -10) node[midway, right, action] {Parity Restore: Transfer max ($1'$) to $B_1$ \\ $B_2$ becomes empty (\fra{case (N1)})};
\node[partition] (P1) at (0, -10.2) {$B_1=\{\mathbf{1}, \mathbf{1'}\}$};

\node[history] (H2) at (7.9, -8.5) {\begin{tabular}{c} \textbf{NODE} $v_2$ \\ \hline \textbf{NODE} $v_3$ \\ \hline \textbf{LOOP} at $v_2$ \\ \hline \textbf{EDGE} $v_2 \to v_3$ \end{tabular}};
\draw[->, dashed, text=gray] (4, -8.5) -- (H2) node[midway, above, font=\tiny, text=gray] {Push};

\draw[->, thick] (P1) -- (0, -12.0) node[midway, right, action] {Remove $\mathbf{1, 1'}$ \\ $B_1$ becomes empty (\fra{case (N0)})};
\node (P0) at (0, -12.5) {\Large $\emptyset$};

\node[history] (H1) at (7.9, -11.5) {\begin{tabular}{c} \textbf{NODE} $v_1$ \\ \hline \textbf{NODE} $v_2$ \\ \hline \textbf{NODE} $v_3$ \\ \hline \textbf{LOOP} at $v_2$ \\ \hline \textbf{EDGE} $v_2 \to v_3$ \end{tabular}};
\draw[->, dashed, text=gray] (4.8, -11.5) -- (H1) node[midway, above, font=\tiny, text=gray] {Push};

\draw[rounded corners, thick, draw=black, fill=gray!5] (-4, -13) rectangle (10, -16.5);
\node[font=\bfseries] at (3, -13.5) {Reconstructed DIG-model Graph};

\node[graphnode] (v1) at (-1, -15) {$v_1$};
\node[below=0.1cm of v1, font=\footnotesize, text=gray] {Step 1 (NODE)};

\node[graphnode] (v2) at (3, -15) {$v_2$};
\node[below=0.1cm of v2, font=\footnotesize, text=gray] {Step 2 (NODE)};

\node[graphnode] (v3) at (7, -15) {$v_3$};
\node[below=0.1cm of v3, font=\footnotesize, text=gray] {Step 3 (NODE)};

\draw[graphedge] (v2) to[out=135,in=45,loop,looseness=5] node[above, font=\footnotesize, text=gray] {Step 4 (LOOP)} (v2);
\draw[graphedge] (v2) -- (v3) node[midway, above, font=\footnotesize, text=gray] {Step 5 (EDGE)};

\end{tikzpicture}
    }
    \caption{Step-by-step execution of \Cref{algo:1,algo:2} on the partition in even subsets $\{\{1, 2'\}, \{1', 2, 4, 5\},\{3, 3', 4', 5'\}\}$, with the order $1\prec1'\prec2\prec2'\prec\cdots$. The blocks are eliminated in the order $B_3,B_2,B_1$, so that $B_i$ corresponds to vertex $v_i$.}
    \label{fig:partition_to_graph}
\end{figure}
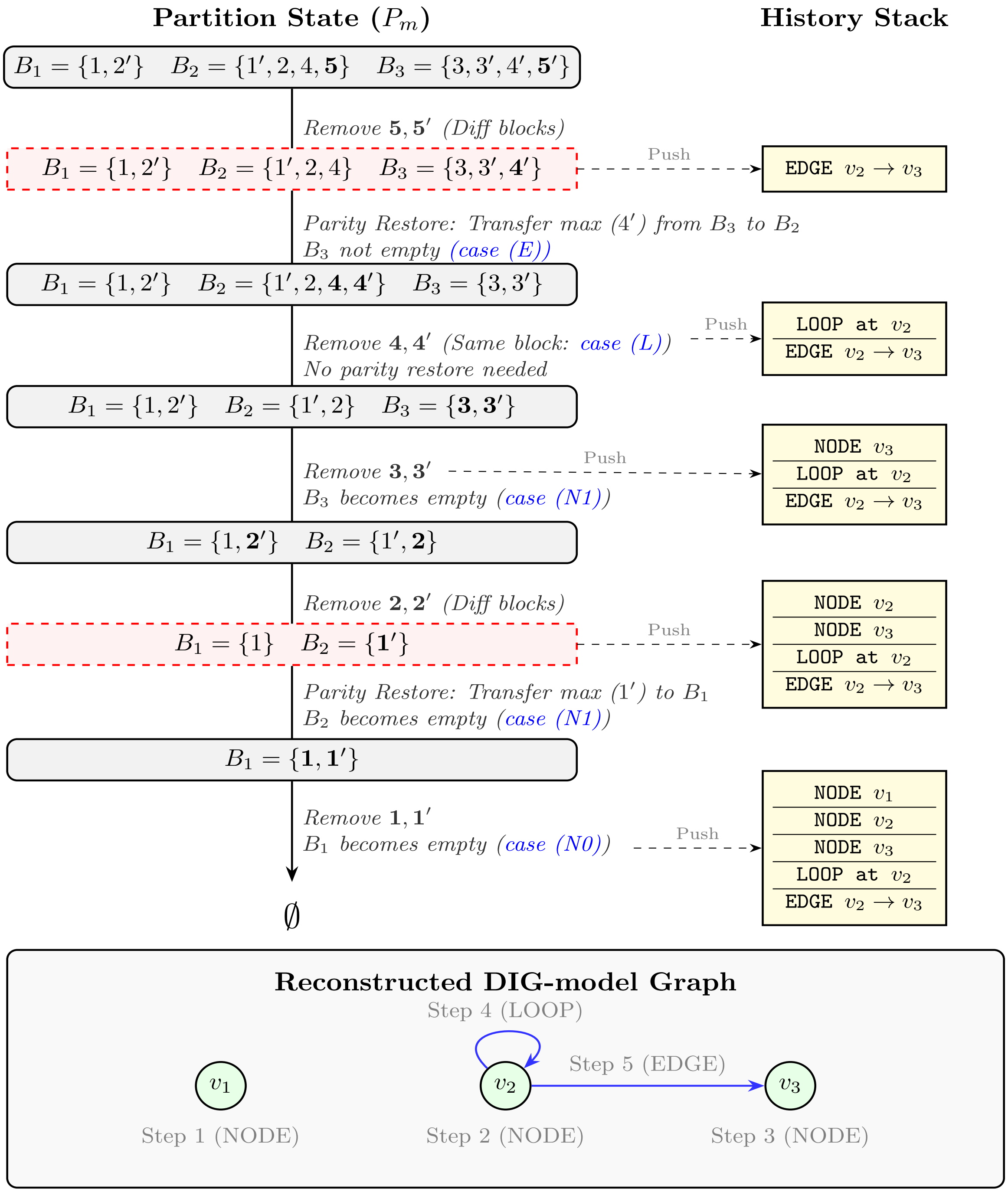

We start with two algorithms working with a data structure called \texttt{History}. 
\Cref{algo:1} extracts a history of graph-building instructions from a partition while applying the insertion operations (N0), (N1), (L), (E), which we introduced in Section~\ref{subsec:directed-even-partitions}, in reverse order and outputs the instructions to \texttt{History}. 
\Cref{algo:2} assembles the contents of \texttt{History} into a DIG-model graph.
See Figure~\ref{fig:partition_to_graph} for the a step-by-step execution of both algorithms on an example.

This two-phased procedure exposes the origin of the many-to-one mapping.
During deconstruction, every time a block is eliminated (a node-defining step), the 
information whether the step was the reverse of (N0) or of (N1) is discarded. Similarly, we discard information on the specific block and chosen element in the cases where we revert (L) or (E). As a graph with $k$ nodes requires exactly $k$ such node-defining steps, offering  $1$, $3$, $\dots$, $2k-1$ choices for the first, second, \dots, $k$-th step, respectively, there are exactly $(2k-1)!!$ distinct pairings that are forgotten during parsing. 
Conversely, the locations of the pairs during the $k^2$ edge and loop steps are preserved within the history, allowing the graph to be deterministically reconstructed. 
Therefore, any single reconstructed graph serves as the representative for exactly $(2k-1)!!$ distinct even-sized partitions.

\begin{algorithm}
\caption{DeconstructPartition($\mathcal{P}$)}\label{algo:1}
\begin{algorithmic}[1]
\REQUIRE A partition $\mathcal{P}$ of $\{1, 1', \dots, n, n'\}$ into $k$ blocks of even size; the elements are ordered by $1\prec1'\prec2\prec\cdots\prec n\prec n'$.
\STATE Initialize \texttt{History} as an empty stack and $L$ as an empty list.
\STATE \COMMENT{Blocks keep their time stamp when elements are removed or transferred.}
\FOR{$m = n$ \textbf{down to} $1$}
    \STATE Let $C$ and $C'$ be the blocks containing $m$ and $m'$, respectively.
    \STATE Remove $m$ and $m'$ from $\mathcal{P}$.
    \IF{$C=C'$}
        \IF{$C$ is empty}
            \STATE \texttt{Instruction} $\leftarrow$ \texttt{NODE}$(C)$; delete $C$ and append $C$ to $L$. \COMMENT{\fra{Case (N0)}}
        \ELSE
            \STATE \texttt{Instruction} $\leftarrow$ \texttt{LOOP}$(C)$. \COMMENT{\fra{Case (L)}}
        \ENDIF
    \ELSE
        \STATE $x\leftarrow\max_\prec(C\cup C')$; let $D\in\{C,C'\}$ be the block containing $x$ and $\bar D$ the other one.
        \STATE Transfer $x$ from $D$ to $\bar D$. \COMMENT{parity-restoring operation}
        \IF{$D$ is empty}
            \STATE \texttt{Instruction} $\leftarrow$ \texttt{NODE}$(D)$; delete $D$ and append $D$ to $L$. \COMMENT{\fra{Case (N1)}}
        \ELSE
            \STATE \texttt{Instruction} $\leftarrow$ \texttt{EDGE}$(C\to C')$. \COMMENT{\fra{Case (E)}}
        \ENDIF
    \ENDIF
    \STATE Push \texttt{Instruction} to \texttt{History}.
\ENDFOR
\STATE The $i$-th block of $L$ receives the vertex label $k+1-i$; replace every block occurring in \texttt{History} by its label.
\RETURN \texttt{History}
\end{algorithmic}
\end{algorithm}

\begin{algorithm}
\caption{BuildGraph(\texttt{History})}
\begin{algorithmic}[1]\label{algo:2}
\REQUIRE The parsed \texttt{History} stack from DeconstructPartition.
\STATE Initialize directed graph $G$ with $0$ vertices.
\STATE Initialize step counter $m \leftarrow 1$.
\WHILE{\texttt{History} is not empty}
    \STATE Pop \texttt{Instruction}.
    \IF{\texttt{Instruction} is \texttt{NODE}$(i)$}
        \STATE Create a new vertex $v_i$ in $G$ with time stamp $m$. \COMMENT{\fra{Cases (N0) and (N1)}}
    \ELSIF{\texttt{Instruction} is \texttt{LOOP}$(i)$}
        \STATE Add a directed loop with time stamp $m$ at vertex $v_i$. \COMMENT{\fra{Case (L)}}
    \ELSIF{\texttt{Instruction} is \texttt{EDGE}$(i\to j)$}
        \STATE Add a directed edge with time stamp $m$ from $v_i$ to $v_j$. \COMMENT{\fra{Case (E)}}
    \ENDIF
    \STATE $m \leftarrow m + 1$.
\ENDWHILE
\RETURN $G$
\end{algorithmic}
\end{algorithm}

\medskip
Next, we turn to the needed bias. Note that the standard exponential generating function for set partitions into even-sized blocks is $\exp(g(t))$, where the inner function $g(t) = \cosh(\sqrt{2t}) - 1$ is associated with the combinatorial class of blocks of even size.   
To inject the biasing factor $1/(2k-1)!!$, we replace the outer exponential series with the function $\cE(y) = \cosh(\sqrt{2y})$, derived from the rescaled Borel transform. 

The series expansion of $\cE(g(t))$ assigns the weight $g(t)^k/((2k-1)!!k!)$ to configurations containing $k$ blocks. 
Based on this algebraic structure, we design a two-stage sampling procedure for any positive real parameter $t$, which can be chosen freely. In practice, it is set to tune the expected number $n$ of steps used to generate the graph. 
First, we draw the integer $k$ from the discrete probability distribution proportional to the Taylor coefficients of $\cE(y)$ evaluated at $y=g(t)$. 
Second, we determine the specific sizes of these $k$ blocks. This is achieved by drawing $k$ independent random variables from the distribution proportional to the coefficients of $g(t)$. 
By constructing a partition from these specific biased distributions and applying \Cref{algo:1,algo:2} afterwards, the resulting distribution over histories of the DIG-model  generated in a given number of steps becomes uniform.

\begin{algorithm}
\caption{Uniform Sampler for DIG-Model Graphs}
\label{alg:uniform_sampler}
\begin{algorithmic}[1]
\REQUIRE A parameter $t > 0$.
\STATE \COMMENT{Phase 1: Biased sampling of the number of blocks $K$}
\STATE Draw $K$ from a discrete distribution proportional to the coefficients of $\cE(y) = \cosh(\sqrt{2y})$ evaluated at $g(t)$: $$\mathbb{P}_t(K = k) = \frac{2^k}{(2k)!} \frac{g(t)^k}{\overrightarrow{B}(t, 1)},\quad \mathrm{for\ }k \ge 0.$$
\STATE \COMMENT{Phase 2: Sampling the sizes of the $K$ even blocks}
\FOR{$i = 1$ \textbf{to} $K$}
    \STATE Draw the size $S_i$ of block $i$ from a distribution proportional to the coefficients of $g(t) = \cosh(\sqrt{2t}) - 1$: $$\mathbb{P}_t(S_i = s_i) = \frac{2^{s_i/2} t^{s_i/2}}{s_i! \,g(t)},\quad \mathrm{for\ even\ }s_i \ge 2.$$
\ENDFOR
\STATE $2n \leftarrow \sum_{i=1}^K S_i$
\STATE \COMMENT{Phase 3: Partition construction and graph projection}
\STATE Draw a uniformly random permutation of $\{1,1',\ldots,n,n'\}$ and cut it into consecutive blocks of sizes $S_1, \dots, S_K$; let $\mathcal{P}$ be the resulting partition.
\STATE $G \leftarrow \text{BuildGraph}(\text{DeconstructPartition}(\mathcal{P}))$
\RETURN $G$
\end{algorithmic}
\end{algorithm}

\begin{theorem}
For any $t > 0$, Algorithm \ref{alg:uniform_sampler} generates a random graph from the DIG-model graph.
Conditioned on the generated graph having $n$ time steps,
the resulting probability distribution is uniform.
\end{theorem}

\begin{proof}
Let $G$ be a history of the DIG-model that is generated in $n$ time steps and has $k$ vertices.
Let $\Pi(G)$ 
denote the equivalence class of even-sized set partitions that are deterministically mapped to $G$ by \Cref{algo:1,algo:2}. 
By the correspondence established in Section~\ref{subsec:directed-even-partitions}, its cardinality is $|\Pi(G)| = (2k-1)!!$.

Fix a partition $\mathcal P$ of $\{1,1',\ldots,n,n'\}$ into $k$ even blocks of sizes $s_1,\ldots,s_k$ (in some order). The algorithm produces $\mathcal P$ if and only if $K=k$, the sequence $(S_1,\ldots,S_k)$ is the sequence of the block sizes of $\mathcal P$ in one of the $k!$ orders of its blocks, and Phase~3 decomposes the $2n$-element ground set into these blocks in this order. 
Since Phase~3 constructs the partition by cutting a uniformly random permutation of length $2n$, generating exactly these specific blocks in this fixed order occurs if and only if the permutation groups the elements of each block together. There are exactly $\prod_i s_i!$ such favorable permutations out of the $(2n)!$ possible ones, so this last event has probability $\lp\prod_i s_i!\rp/(2n)!$. Since $\sum_i s_i/2=n$, we obtain
\begin{align*}
\mathbb{P}_t(\mathcal{P})
&= k!\cdot\frac{2^k g(t)^k}{(2k)!\overrightarrow{B}(t,1)}\cdot\frac{(2t)^n}{g(t)^k\prod_i s_i!}\cdot\frac{\prod_i s_i!}{(2n)!}
= \frac{1}{\overrightarrow{B}(t,1)} \frac{(2t)^n}{(2n)!} \frac{1}{(2k-1)!!},
\end{align*}
because $k!2^k/(2k)!=1/(2k-1)!!$.

The overall probability of sampling the final graph $G$ is the sum of the probabilities of generating any partition that is mapped to $G$:
\begin{equation}
    \mathbb{P}_t(G) = \sum_{\mathcal{P} \in \Pi(G)} \mathbb{P}_t(\mathcal{P})
    = (2k-1)!! \left( \frac{1}{\overrightarrow{B}(t,1)} \frac{(2t)^n}{(2n)!} \frac{1}{(2k-1)!!} \right)
    = \frac{1}{\overrightarrow{B}(t,1)} \frac{(2t)^n}{(2n)!}.
\end{equation}

This probability $\mathbb{P}_t(G)$ depends exclusively on $n$ and the parameter $t$. It is entirely independent of the number of nodes $k$. Consequently, conditioning this generation process on the output being generated in $n$ time steps guarantees a uniform distribution over all the DIG-model graphs with $n$ time steps.
\end{proof}

\begin{remark}
Summing over all histories of length $n$ gives
\[
\mathbb P_t(N=n)=\frac{c_n(2t)^n}{(2n)!\overrightarrow{B}(t,1)}
=\frac{t^n[t^n]\overrightarrow{B}(t,1)}{\overrightarrow{B}(t,1)}.
\]
To sample a history of exact length $n$ by rejection, it is natural to choose $t=\rho_n$, the solution of $a(\rho_n)=n$ from the proof of Theorem~\ref{directed-asymptotic}; then $\EE_t[N]=a(t)=n$. By Hayman's theorem and \eqref{directed-b-at-saddle},
\[
\mathbb P_{\rho_n}(N=n)\sim\frac{1}{\sqrt{2\pi b(\rho_n)}}\sim\frac{1}{\sqrt{\pi n(w_n+1)}},
\]
so that the expected number of trials is of order $\sqrt{n\log n}$.
\end{remark}

\section{Conclusion}

In this work we have proposed two random graph growth models where the growth is achieved by adding a new vertex or a new edge uniformly at random. With the UIG-model and the DIG-model, two natural models have been chosen for our investigation as prototypes for an undirected and a directed model. The two models are interesting in their own right, as they exhibit the need for a fairly different mathematical treatment. In both models, the number of objects of size $n$ (where size is the length of the history) grows so rapidly that the generating functions that we encountered have a radius of convergence equal to zero. In the UIG-model, the expression as a sum together with the discrete Laplace method was the clue. The DIG-model required a rescaled Borel transform to make Hayman's theory of saddle-point asymptotics and some of its extensions work. In both cases we see asymptotic expressions involving the Lambert W function. 

Of course, two other obvious models are undirected graphs with loops allowed and directed graphs without loops. But we expect that one of the two approaches will work in either case. Other models in the framework of independent growth of the number of vertices and the number of edges can be obtained by imposing certain constraints of our model. 

The most natural parameter in these models is of course the number of vertices of the graph generated by a history, for which we show a central limit theorem in both models. Other parameters are conceivable and in the extended abstract of this work~\cite{AzzouzBodiniDurandGittenberger2026} we treated some of them briefly. It is also worth investigating whether the correspondence between the DIG-model and set partitions can lead to interesting asymptotic results, either for the DIG-model graphs or for partitions. A natural line of future research would be the question whether a modification of the rules to generate a history, for instance looking at hypergraphs instead of graphs, will lead to structures that can be analyzed in a similar way and lead to similar asymptotics. Maybe such structures satisfy a correspondence to set partitions with different divisibility properties for their block sizes. 

A brief discussion was done when introducing weights on vertex and edge creation. Unfortunately, the approach seems not to be well-suited for problems of that kind. So, it remains an open problem how to introduce weights in a proper way to allow for tuning the proportion of the number of edges and vertices in a history, such that the resulting graph class is amenable to a complete asymptotic analysis, including the description of possible phase transitions. 

As the transformed generating function is nearly bivariate Hayman admissible it is worthwhile to investigate this further, In~\cite{AzzouzBodiniDurandGittenberger2026} we customized bivariate Hayman admissibility to our needs and were able to carry over the main theorem to our case. This gives a full bivariate normal limit law for the joint distribution of the number of vertices and the length of the history. A question is whether the slight extension of the concept developed in~\cite{GittenbergerMandlburger2006} mentioned above or a similar one allows for inclusion of further parameters, yielding multivariate normal limit laws based on multivariate Hayman admissibility. 

Finally, we remark that this work serves as another example demonstrating that 
encountering combinatorial structures with rapidly growing counting sequences — whose generating functions are merely formal power series — is not the end of the story in analytic combinatorics. There are a lot of tools available to rigorously treat divergent series in a proper way to squeeze out precise asymptotic information from them.

\printbibliography

@article{GittenbergerMandlburger2006,
  author  = {Bernhard Gittenberger and Johannes Mandlburger},
  title   = {Hayman admissible functions in several variables},
  journal = {The Electronic Journal of Combinatorics},
  volume  = {13},
  number  = {1},
  pages   = {R106},
  year    = {2006},
  doi     = {10.37236/1132}
}

@article{Hayman1956,
  author  = {Hayman, W. K.},
  title   = {A generalisation of Stirling's formula},
  journal = {Journal f{\"u}r die reine und angewandte Mathematik},
  volume  = {196},
  year    = {1956},
  pages   = {67--95},
  doi     = {10.1515/crll.1956.196.67},
  url     = {https://www.degruyterbrill.com/document/doi/10.1515/crll.1956.196.67/html}
}

@inproceedings{berger2006framework,
  author    = {Berger-Wolf, Tanya Y. and Saia, Jared},
  title     = {A framework for analysis of dynamic social networks},
  booktitle = {Proceedings of the 12th ACM SIGKDD International Conference on Knowledge Discovery and Data Mining (KDD '06)},
  pages     = {523--528},
  publisher = {ACM},
  year      = {2006},
  doi       = {10.1145/1150402.1150462}
}

@inproceedings{behr2020rewriting,
  author    = {Behr, Nicolas and Krivine, Jean},
  title     = {Rewriting theory for the life sciences: a unifying theory of {CTMC} semantics},
  booktitle = {Graph Transformation (ICGT 2020)},
  series    = {Lecture Notes in Computer Science},
  volume    = {12150},
  pages     = {185--202},
  publisher = {Springer},
  year      = {2020},
  doi       = {10.1007/978-3-030-51372-6_11}
}

@article {CooperFrieze03,
    AUTHOR = {Cooper, Colin and Frieze, Alan},
     TITLE = {A general model of web graphs},
   JOURNAL = {Random Structures Algorithms},
  FJOURNAL = {Random Structures \& Algorithms},
    VOLUME = {22},
      YEAR = {2003},
    NUMBER = {3},
     PAGES = {311--335},
      ISSN = {1042-9832,1098-2418},
   MRCLASS = {05C80 (05C07 60C05 68R10)},
  MRNUMBER = {1966545},
MRREVIEWER = {Hsien-Kuei\ Hwang},
       DOI = {10.1002/rsa.10084},
       URL = {https://doi.org/10.1002/rsa.10084},
}

@book {Durrett07,
    AUTHOR = {Durrett, Rick},
     TITLE = {Random {G}raph {D}ynamics},
    SERIES = {Cambridge Series in Statistical and Probabilistic Mathematics},
    VOLUME = {20},
 PUBLISHER = {Cambridge University Press, Cambridge},
      YEAR = {2007},
     PAGES = {x+212},
      ISBN = {978-0-521-86656-9},
   MRCLASS = {05C80 (05-02 60-02 60C05 60G50 60K35 82C41)},
  MRNUMBER = {2271734},
MRREVIEWER = {Michael\ Krivelevich},
}

@article{farber2024random,
  author  = {Farber, Michael and Gnedin, Alexander and Mannan, Wajid},
  title   = {A random graph growth model},
  journal = {Bulletin of the London Mathematical Society},
  volume  = {56},
  number  = {2},
  pages   = {662--680},
  year    = {2024},
  doi     = {10.1112/blms.12957}
}

@book {ChungLu06,
    AUTHOR = {Chung, Fan and Lu, Linyuan},
     TITLE = {Complex {G}raphs and {N}etworks},
    SERIES = {CBMS Regional Conference Series in Mathematics},
    VOLUME = {107},
 PUBLISHER = {Published for the Conference Board of the Mathematical
              Sciences, Washington, DC; by the American Mathematical
              Society, Providence, RI},
      YEAR = {2006},
     PAGES = {viii+264},
      ISBN = {978-0-8218-3657-6},
   MRCLASS = {05C85 (68R10 94C15)},
  MRNUMBER = {2248695},
MRREVIEWER = {Colin\ D.\ Cooper},
       DOI = {10.1090/cbms/107},
       URL = {https://doi.org/10.1090/cbms/107},
}

@book {Bonato08,
    AUTHOR = {Bonato, Anthony},
     TITLE = {A {C}ourse on the {W}eb {G}raph},
    SERIES = {Graduate Studies in Mathematics},
    VOLUME = {89},
 PUBLISHER = {American Mathematical Society, Providence, RI; Atlantic
              Association for Research in the Mathematical Sciences (AARMS),
              Halifax, NS},
      YEAR = {2008},
     PAGES = {xii+184},
      ISBN = {978-0-8218-4467-0},
   MRCLASS = {05-01 (05C80 05C90 68R10 94C15)},
  MRNUMBER = {2389013},
MRREVIEWER = {Abraham\ D.\ Flaxman},
       DOI = {10.1090/gsm/089},
       URL = {https://doi.org/10.1090/gsm/089},
}

@article{blasiak2010combinatorial,
    AUTHOR = {Blasiak, Pawel and Flajolet, Philippe},
     TITLE = {Combinatorial models of creation-annihilation},
   JOURNAL = {S\'{e}m. Lothar. Combin.},
  FJOURNAL = {S\'{e}minaire Lotharingien de Combinatoire},
    VOLUME = {65},
      YEAR = {2010},
     PAGES = {Art. B65c, 78},
      ISSN = {1286-4889}
}

@article{andrews2011legendre,
  title={The {L}egendre--{S}tirling numbers},
  author={Andrews, George E and Gawronski, Wolfgang and Littlejohn, Lance L},
  journal={Discrete Mathematics},
  volume={311},
  number={14},
  pages={1255--1272},
  year={2011},
  publisher={Elsevier}
}

@article{andrews2009combinatorial,
  title={A combinatorial interpretation of the {L}egendre--{S}tirling numbers},
  author={Andrews, George and Littlejohn, Lance},
  journal={Proceedings of the American Mathematical Society},
  volume={137},
  number={8},
  pages={2581--2590},
  year={2009}
}

@book{comtet2012advanced,
  title      = {Advanced Combinatorics: The Art of Finite and Infinite Expansions},
  author     = {Comtet, Louis},
  translator = {Nienhuys, J. W.},
  year       = {2012},
  publisher  = {Springer Science \& Business Media},
  address    = {Dordrecht},
  doi        = {10.1007/978-94-010-2196-8},
  note       = {Reprint of the original 1974 edition}
}

@article{Duchon2004,
  title={Boltzmann samplers for the random generation of combinatorial structures},
  author={Duchon, Philippe and Flajolet, Philippe and Louchard, Guy and Schaeffer, Gilles},
  journal={Comb. Probab. Comput.},
  volume={13},
  number={4-5},
  pages={577--625},
  year={2004},
  publisher={Cambridge University Press}
}

@article{carlitz1963divided,
  title={The divided central differences of zero},
  author={Carlitz, L and Riordan, John},
  journal={Canadian Journal of Mathematics},
  volume={15},
  pages={94--100},
  year={1963},
  publisher={Cambridge University Press}
}

@article{everitt2002legendre,
  title={Legendre polynomials, {L}egendre--{S}tirling numbers, and the left-definite spectral analysis of the Legendre differential expression},
  author={Everitt, WN and Littlejohn, LL and Wellman, R},
  journal={Journal of Computational and Applied Mathematics},
  volume={148},
  number={1},
  pages={213--238},
  year={2002},
  publisher={Elsevier}
}

@article{erdHos1960evolution,
  author  = {Erd{\H{o}}s, Paul and R{\'e}nyi, Alfr{\'e}d},
  title   = {On the evolution of random graphs},
  journal = {Publication of the Mathematical Institute of the Hungarian Academy of Sciences},
  volume  = {5},
  pages   = {17--61},
  year    = {1960}
}

@article{albert2002statistical,
  title={Statistical mechanics of complex networks},
  author={Albert, R{\'e}ka and Barab{\'a}si, Albert-L{\'a}szl{\'o}},
  journal={Reviews of Modern Physics},
  volume={74},
  number={1},
  pages={p. 47},
  year={2002},
  publisher={APS}
}

@book{flajolet2009analytic,
  title={Analytic {C}ombinatorics},
  author={Flajolet, Philippe and Sedgewick, Robert},
  year={2009},
  publisher={Cambridge University Press}
}

@book {JLR00,
    AUTHOR = {Janson, Svante and {\L}uczak, Tomasz and Rucinski, Andrzej},
     TITLE = {Random graphs},
    SERIES = {Wiley Interscience Series in Discrete Mathematics and
              Optimization},
 PUBLISHER = {Wiley-Interscience, New York},
      YEAR = {2000},
     PAGES = {xii+333},
      ISBN = {0-471-17541-2},
   MRCLASS = {05C80 (60C05 82B41)},
  MRNUMBER = {1782847},
MRREVIEWER = {Mark\ R.\ Jerrum},
       DOI = {10.1002/9781118032718},
       URL = {https://doi.org/10.1002/9781118032718},
}

@book {Bollobas01,
    AUTHOR = {Bollob\'{a}s, B\'{e}la},
     TITLE = {Random Graphs},
    SERIES = {Cambridge Studies in Advanced Mathematics},
    VOLUME = {73},
   EDITION = {Second},
 PUBLISHER = {Cambridge University Press, Cambridge},
      YEAR = {2001},
     PAGES = {xviii+498},
      ISBN = {0-521-80920-7},
      NOTE = {Paperback ISBN: 0-521-79722-5},
   MRCLASS = {05C80 (60C05)},
  MRNUMBER = {1864966},
       DOI = {10.1017/CBO9780511814068},
       URL = {https://doi.org/10.1017/CBO9780511814068},
}

@book{Ovidiu2008,
  author    = {Ovidiu Costin},
  title     = {Asymptotics and Borel Summability},
  publisher = {CRC Press / Taylor \& Francis},
  year      = {2008},
  series    = {Monographs and Surveys in Pure and Applied Mathematics},
  isbn      = {978-1420070316}
}

@article{KimKimCentralBell,
  author  = {Kim, Taekyun and Kim, Dae San},
  title   = {A note on central Bell numbers and polynomials},
  journal = {Russian Journal of Mathematical Physics},
  volume  = {27},
  number  = {1},
  pages   = {76--81},
  year    = {2020},
  doi     = {10.1134/S1061920820010070}
}

@book{Balser1994,
  author    = {Balser, Werner},
  title     = {From Divergent Power Series to Analytic Functions:
               Theory and Application of Multisummable Power Series},
  series    = {Lecture Notes in Mathematics},
  volume    = {1582},
  publisher = {Springer},
  address   = {Berlin, Heidelberg},
  year      = {1994},
  doi       = {10.1007/BFb0073564},
  isbn      = {978-3-540-58268-7}
}

@article{MichalikTkacz2019,
  author  = {Michalik, S{\l}awomir and Tkacz, Bo{\.z}ena},
  title   = {The Stokes phenomenon for some moment partial differential equations},
  journal = {Journal of Dynamical and Control Systems},
  volume  = {25},
  pages   = {573--598},
  year    = {2019},
  doi     = {10.1007/s10883-018-9424-9}
}

@misc{NIST:DLMF,
  key          = "{\relax DLMF}",
  title        = "{\it NIST Digital Library of Mathematical Functions}",
  howpublished = "\url{https://dlmf.nist.gov/}, Release 1.2.8 of 2026-09-15",
  url          = "https://dlmf.nist.gov/",
  note         = "F.~W.~J. Olver, A.~B. {Olde Daalhuis}, D.~W. Lozier,
                  B.~I. Schneider, R.~F. Boisvert, C.~W. Clark,
                  B.~R. Miller, B.~V. Saunders, H.~S. Cohl,
                  and M.~A. McClain, eds."
}

@article{DrmotaGittenbergerKlausner2005,
  author  = {Drmota, Michael and Gittenberger, Bernhard and Klausner, Thomas},
  title   = {Extended admissible functions and Gaussian limiting distributions},
  journal = {Mathematics of Computation},
  volume  = {74},
  number  = {252},
  pages   = {1953--1966},
  year    = {2005},
  doi     = {10.1090/S0025-5718-05-01744-8},
  url     = {https://doi.org/10.1090/S0025-5718-05-01744-8}
}

@inproceedings{AzzouzBodiniDurandGittenberger2026,
    AUTHOR = {Azzouz, Nadja and Bodini, Olivier and Durand, Francis and
              Gittenberger, Bernhard},
     TITLE = {Asymptotic analysis of generating functions arising from
              dynamic graphs},
 BOOKTITLE = {37th {I}nternational {C}onference on {P}robabilistic,
              {C}ombinatorial and {A}symptotic {M}ethods for the {A}nalysis
              of {A}lgorithms},
    SERIES = {LIPIcs. Leibniz Int. Proc. Inform.},
    VOLUME = {381},
     PAGES = {Paper No. 30, 16},
 PUBLISHER = {Schloss Dagstuhl. Leibniz-Zent. Inform., Wadern},
      YEAR = {2026},
      ISBN = {978-3-95977-435-2},
%  author    = {Azzouz, Nadja and Bodini, Olivier and Durand, Francis and Gittenberger, Bernhard},
%  title     = {Asymptotic analysis of generating functions arising from dynamic graphs},
%  booktitle = {37th International Conference on Probabilistic, Combinatorial and Asymptotic Methods for the Analysis of Algorithms (AofA 2026)},
%  year      = {2026},
%  address   = {Munich, Germany},
  doi       = {10.4230/LIPIcs.AofA.2026.30},
%  note      = {HAL: hal-05727365}
}

@article{corless1996lambertw,
  author  = {Corless, Robert M. and Gonnet, Gaston H. and Hare, David E. G. and Jeffrey, David J. and Knuth, Donald E.},
  title   = {On the Lambert W function},
  journal = {Advances in Computational Mathematics},
  volume  = {5},
  number  = {4},
  pages   = {329--359},
  year    = {1996},
  doi     = {10.1007/BF02124750}
}

@article {klazar2003,
    AUTHOR = {Klazar, Martin},
     TITLE = {Bell numbers, their relatives, and algebraic differential
              equations},
   JOURNAL = {J. Combin. Theory Ser. A},
  FJOURNAL = {Journal of Combinatorial Theory. Series A},
    VOLUME = {102},
      YEAR = {2003},
    NUMBER = {1},
     PAGES = {63--87},
      ISSN = {0097-3165,1096-0899},
       DOI = {10.1016/S0097-3165(03)00014-1},
       URL = {https://doi.org/10.1016/S0097-3165(03)00014-1},
}
\end{document}